\documentclass[12pt]{article}
\usepackage{amsmath,amsfonts}

\newtheorem{thm}{Theorem}
\newtheorem{lemma}{Lemma}

\newcommand{\pfover}{\ \ $\diamond$}

\newcommand{\subsub}[3]{
\displaystyle\mathop{\displaystyle #1_{#2}}_{#3}
}

\newcommand{\CC}{{\mathbb C}}
\newcommand{\kk}{{\mathbf k}}
\newcommand{\PP}{{\mathbb P}}

\newcommand{\bb}{{\mathbf b}}
\newcommand{\cc}{{\mathbf c}}
\newcommand{\Flag}{\mathop{\rm Flag}}
\newcommand{\Ker}{\mathop{\rm Ker}}
\newcommand{\GL}{\mathop{\rm GL}}
\newcommand{\GLn}{{\GL_n}}
\newcommand{\sh}[1]{\!#1\!}
\newcommand{\bdot}{\bullet}
\renewcommand{\leadsto}{\ \Longrightarrow\ }
\newcommand{\C}[1]{{\bigcirc\hspace{-.75em}#1\,}}
\renewcommand{\d}{{\cdot}}
\renewcommand{\max}{{\mathop{\rm max}{}}}
\renewcommand{\min}{{\mathop{\rm min}{}}}
\newcommand{\id}{\text{\rm id}}

\newcommand{\smallstack}[2]
{\,\mbox{\small $\stackrel{\rm #1}{#2}$}\,}
\newcommand{\ldeg}{\smallstack{deg}{\leq}}

\newcommand{\lrk}{\smallstack{rk}{\leq}}
\newcommand{\lmv}{\smallstack{mv}{\leq}}

\newcommand{\sldeg}{\smallstack{deg}{<}}

\newcommand{\slrk}{\smallstack{rk}{<}}

\renewcommand{\lhd}{<\!\!\!\!\cdot\ }
\newcommand{\move}[1]{ \stackrel{\rm #1}{\lhd} }

\newcommand{\tM}{\tilde{M}}

\newcommand{\tD}{\tilde{\Delta}}
\newcommand{\tm}{\tilde{m}}
\newcommand{\tr}{\tilde{r}}
\newcommand{\D}{\Delta}

\newcommand{\rb}[1]{r_{\langle#1\rangle}}
\newcommand{\trb}[1]{\tr_{\langle#1\rangle}}

\newcommand{\bR}{\bar R}

\newcommand{\cbullet}
 { \mbox{\! $\bigcirc\hspace{-.75em}\bullet $ } }
\renewcommand{\d}{\ \bullet\ }
\renewcommand{\c}{\cbullet}
\newcommand{\n}{{}}

\newcommand{\m}[4]
{
\begin{array}{|@{\!\, }c@{\!}c@{\!\, }|}
\hline
#1\, &\, #2\\
#3\, &\, #4\\
\hline
\end{array}
}
\newcommand{\M}[9]
{
\begin{array}{|@{\!\, }c@{\!}c@{\!}c@{\!\, }|}
\hline
#1&\ #2\ &#3\\
#4&\ #5\ &#6\\
#7&\ #8\ &#9\\
\hline
\end{array}
}

\begin{document}

\centerline{\Large\bf BRUHAT ORDER FOR}\vspace{.2em}
\centerline{\Large\bf TWO FLAGS AND A LINE}
\vspace{1em}
\centerline{Peter Magyar}
\vspace{1em}
\centerline{September 1, 2002}

\begin{abstract} \noindent
The classical Ehresmann-Bruhat order
describes the possible degenerations of a pair of flags
in a linear space $V$ under linear transformations of $V$; 
or equivalently, it describes the closure of 
an orbit of $\GL(V)$ acting diagonally on 
the product of two flag varieties.  

We consider the degenerations of
a triple consisting of two flags and a line,
or equivalently the closure of an orbit
of $\GL(V)$ acting diagonally on the product
of two flag varieties and a projective space.
We give a simple rank criterion to decide whether
one triple can degenerate to another.
We also classify the minimal degenerations,
which involve not only reflections 
(i.e., transpositions) in the 
Weyl group $S_n$,\ $n\sh=\dim(V)$, 
but also cycles of arbitrary length.  
Our proofs use only elementary
linear algebra and combinatorics
\end{abstract}

\section{Introduction}

\subsection{A Line and two flags}

We shall deal with certain
configurations of linear subspaces 
in $\CC^n$ (or any vector space $V$).
A configuration $F=(A,B_\bdot,C_\bdot)$
consists of a line $A\subset\CC^n$
and two flags of subspaces of 
fixed dimensions, 
$B_\bdot = (B_1\subset B_2\subset
\cdots\subset\CC^n)$ and
$C_\bdot = (C_1\subset C_2\subset
\cdots\subset\CC^n)$.  In this Introduction,
we restrict ourselves to the case in which
$B_\bdot, C_\bdot$ are full flags:
$\dim B_i=\dim C_i = i$ for $i=0,1,2,\ldots,n$.

Our aim is to describe such configurations 
up to a linear change of coordinates in $\CC^n$,
and the ways in which more generic 
configurations can degenerate to more
special ones.  One could ask this question
for configurations of arbitrarily many flags;
however in general it is `wild' problem.
The distinguishing feature of our case is
that there are only {\it finitely
many} configuration types
$F=(A,B_\bdot,C_\bdot)$, as we 
showed in a previous work \cite{mwz1}
with J.~Weyman and 
A.~Zelevinsky.\footnote{This fact was also noted by Brion \cite{brion}.
More generally, our work 
\cite{mwz1} uses the theory of
quiver representations 
to classify {\it all} products of 
partial flag varieties with finitely 
many orbits of $\GLn$.
See also \cite{mwz2} for the case of the symplectic group $\mathop{\rm Sp}_{2n}$.}

For example, there exists a 
{\it most generic} type $F_\max$,
which degenerates to all other types.
It consists of those configurations
which can be written in terms of 
some basis $v_1,\ldots,v_n$ of $\CC^n$ as:
$$
A = \langle v_1+v_2+\cdots+v_n\rangle,
\qquad
B_i=\langle v_1,v_2,\ldots,v_i\rangle,
\qquad
C_i=\langle v_n,v_{n-1},\ldots,v_{n-i+1}\rangle.
$$
(Here $\langle v_1, v_2,\ldots\rangle$
means the linear span of $v_1,v_2,\ldots$.)
There is also a {\it most special}
configuration type $F_\min$:
$$
A = \langle v_1\rangle,
\qquad
B_i=C_i=\langle v_1,v_2,\ldots,v_i\rangle.
$$

The configurations of a more generic type 
can be made to degenerate to more special
ones by letting some of the basis vectors
$v_i$ approach each other, so that in the 
limit some of the spaces $A, B_i, C_j$
increase their 
intersections.\footnote{ 
That is, $F'$ degenerates to $F$ if we can find
a continuous family of configurations $(A(\tau),B_\bdot(\tau),C_\bdot(\tau))$
indexed by a parameter $\tau\in\CC$, 
such that the configurations for 
$\tau\neq 0$ are all of type $F'$, but
for $\tau=0$ we enter type $F$.}   
Geometrically, a configuration type is an 
orbit of $\GLn(\CC)$ acting diagonally on 
the product $\PP^{n-1}\sh\times\Flag(\CC^n)\sh\times
\Flag(\CC^n)$, with $F_\max$ the open orbit and 
$F_\min$ the unique closed orbit.  
Degeneration of configuration 
types means the topological closure 
of a large orbit contains a smaller orbit.

We seek a simple combinatorial description 
of all degenerations.  The trivial
case of $n=2$ is illustrated by a diagram in \S\ref{examples}.

Our problem is directly analogous to 
the classical case 
in which the configurations consist 
of two flags only: $F=(B_\bdot,C_\bdot)$.  
This theory originated with
Schubert and Ehresmann; a good introduction is \cite{fulton}.
In this case, the configurations 
(up to change of basis in $\CC^n$)
 correspond to permutations $w\in S_n$:
the configuration type $F_w$ consists of 
the double flags which can be written as:
$$
B_i=\langle v_1,v_2,\ldots,v_i\rangle,
\qquad
C_i=\langle v_{w(1)},v_{w(2)},\ldots,v_{w(i)}\rangle
$$
for some basis $v_1,\ldots,v_n$ of $\CC^n$.
A configuration type $F_w$ is a degeneration
of another $F_y$ exactly if $w\leq y$
in the {\it Bruhat order} on $S_n$.
Namely, $w\leq y$ iff
$$
\#(\,[i]\sh\cap w[j]\,)\geq \#(\,[i]\sh\cap y[j]\,)
$$
for all $1\leq i,j\leq n$, where
$[i]:=\{1,2,\ldots,i\}$ and 
$w[j]:=\{w(1),w(2),\ldots,w(j)\}$.
This {\it tableau criterion} has the geometric
meaning:
$$
\#(\,[i]\sh\cap w[j]\,)=\dim(B_i\cap C_j)
$$
for $(B_\bdot,C_\bdot)$ of type $F_w$.
The more special configuration $F_w$ has larger
intersections among its spaces than 
the more generic $F_y$.

We can also describe the classical Bruhat order
in terms of its covers:
$w\lhd y$ iff $y=(i_0,i_1)\cdot w$
for some transposition $(i_0,i_1)\in S_n$
and $\ell(y)=1+\ell(w)$, where
$\ell(w)$ is the number of inversions of $w$.
We can picture this definition in terms of the
permutation matrices $M_w=(m_{ij})$
and $M_y=(m'_{ij})$,
where $m_{ij}:=\delta_{w(i),j}$,\
$m'_{ij}:=\delta_{y(i),j}$.
Then $w\lhd y$ means that we have a pair
of entries $m_{i_0j_0}=m_{i_1j_1}=1$
with $(i_0,j_0)$ northwest of $(i_1,j_1)$,
and no other 1's in the rectangle
$[i_0,i_1]\sh\times[j_0,j_1]$;
and we flip these two `diagonal' entries
in $M_w$ to the corresponding
anti-diagonal, obtaining $M_y$:
$$
w=\hspace{-.5em}\begin{array}{c@{\!}l}
&\hspace{1.7em}j_0\hspace{2.3em}j_1\\
\begin{array}{c}   i_0\\[1.0em] i_1\\ \end{array}
&\begin{array}{|c@{\!}ccc@{\!}c|}
\hline
&\vdots&&\vdots&\\[-.3em]
\cdots\ & 1 & \cdots & 0 &\ \cdots\\[-.5em]
&\vdots&0&\vdots&\\
\cdots\ & 0 &\cdots& 1 &\ \cdots\\[-.5em]
&\vdots&&\vdots&\\
\hline
\end{array}
\end{array}
\ \lhd\  
y=\hspace{-.5em}\begin{array}{c@{\!}l}
&\hspace{1.7em}j_0\hspace{2.3em}j_1\\
\begin{array}{c}   i_0\\[1.0em] i_1\\ \end{array}
&\begin{array}{|c@{\!}ccc@{\!}c|}
\hline
&\vdots&&\vdots&\\[-.3em]
\cdots\ & 0 & \cdots & 1 &\ \cdots\\[-.5em]
&\vdots&0&\vdots&\\
\cdots\ & 1 &\cdots& 0 &\ \cdots\\[-.5em]
&\vdots&&\vdots&\\
\hline
\end{array}
\end{array}\,;
$$
or in compact notation, with 1 replaced
by $\bullet$ and all unaffected rows and
columns omitted:
$$
\m{\d}{} {}{\d}\ \lhd\ \m{}{\d} {\d}{}\ .
$$
In terms of transpositions:
$y=(i_0,i_1)\sh\cdot w = w\sh\cdot(j_0,j_1)$.

We give a full exposition and proof
of these classical results
in \S\ref{sec:abruhat1}, \S\ref{sec:abruhat2}.

\subsection{Bruhat order}

Let us return to our case of a line and
two flags.
As we showed in \cite{mwz1},
we can index our configuration types
by {\it decorated permutations} $(w,\Delta)$,
where $\Delta=\{j_1\sh<j_2\sh<\cdots\sh<j_t\}$
is any non-empty descending subsequence of $w$,
meaning $w(j_1)>w(j_2)>\cdots> w(j_t)$.
In the corresponding configuration $F_{w,\Delta}$,
the permutation $w$ describes the relative positions of $B_\bdot$
and $C_\bdot$ in terms of a basis
$v_1,\ldots,v_n$, just as before;
and $\Delta$ defines the extra line:
$$
A=\langle v_{j_1}+v_{j_2}+\cdots+v_{j_t}\rangle\,.
$$

Thus, the generic $F_\max$ is $F_{w,\Delta}$
for $w=w_0=\underline{n,n\sh-1,\ldots,2,1}$,
the longest permutation, and 
$\Delta=\{1,2,\ldots,n\}$.
The most special $F_\min$ is $F_{w,\Delta}$
for $w=\id=\underline{1,2,\ldots,n}$
and $\Delta=\{1\}$.
We can picture a decorated permutation
as a permutation matrix with circles
around the positions $(w(j),j)$ for
$j\in\Delta$.  For example, 
{\tiny $\begin{array}{|@{\!}c@{\!}c@{\!}c@{\!}|}\hline 
 &\ \bullet\ &  \\
 & &\cbullet\\
\cbullet& & \\
\hline\end{array}
$}\ \ 
corresponds to $w=\underline{312}$ , 
$\Delta=\{1,3\}$.

We once again have a degeneration
or Bruhat order, described combinatorially
by a tableau criterion in terms of
certain {\it rank numbers} which measure
intersections of spaces in a configuration
$(A, B_\bdot,C_\bdot)$ in $F_{w,\D}$.
Namely, let
$$
r_{ij}(w):=\dim(B_i\cap C_j)=
\#(\,[i]\sh\cap w[j]\,)
$$
as before, and
$$\begin{array}{rcl}
\rb{ij}(w,\Delta)&:=&\dim(B_i\cap C_j)+
\dim(\,A\cap(B_i\sh+C_j)\,)\\[.3em]
&\,=&\#(\,[i]\sh\cap w[j]\,)+\,
\delta_{ij}(w,\Delta)\,,
\end{array}$$
where
$$
\delta_{ij}(w,\Delta):=\left\{
\begin{array}{cl}
1&\text{if for all}\ k\in\Delta,\ \ 
k\leq i\ \text{or}\ w(k)\leq j\\
0&\text{otherwise}\ .
\end{array}\right.
$$
We can realize this in terms of linear algebra
by defining $\phi_{ij}:B_i\sh\times C_j\to
\CC^n/A$,\ $(v_1,v_2)\mapsto v_1\sh+v_2\mod A$:
then 
$
\rb{ij}(w,\D)=\dim\Ker \phi_{ij}\,.
$
These definitions are suggested by
quiver theory: see \S\ref{structure} below.
We will show that our geometric degeneration order
has the following combinatorial description: 
$$
(w,\Delta)\leq(y,\Gamma)\quad
\Leftrightarrow
\quad \left\{
\begin{array}{c}
r_{ij}(w,\Delta)\geq r_{ij}(y,\Gamma)\\
\rb{ij}(w,\Delta)\geq \rb{ij}(y,\Gamma)\\
\text{for all}\ 0\leq i,j\leq n
\end{array}\right.
\ .
$$ 

Finally, we can classify the covers 
$(w,\Delta)\move{}(y,\Gamma)$ of our
new Bruhat order.  Remarkably,
in many of the cases below 
the pair $w<y$ is {\it not} a cover 
in the classical Bruhat order.
We describe the covers in terms of certain flipping
moves which we write in compact notation
(again, with all unaffected 
rows and columns omitted).
We describe how a more generic 
configuration $(A, B_\bdot,C_\bdot)$ (on the right)
degenerates to a more special configuration
(on the left).
\\[.5em]
{\sc move} (i)\ The line $A$ moves into one of the
spaces $B_i+ C_j$, leaving $B_\bdot$, $C_\bdot$
unchanged:
$$
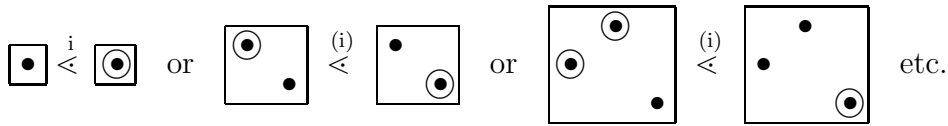
\begin{array}{|@{\!}c@{\!}|}\hline\,\d\,\\ \hline\end{array}
\move{i}
\begin{array}{|@{\!}c@{\!}|}\hline\,\c\,\\ \hline\end{array}
\quad\text{or}\quad
\m{\c}{}{}{\d}\ \move{(i)}\
\m{\d}{}{}{\c}
\quad\text{or}\quad
\M{}{\c}{} {\c}{}{} {}{}{\d}
\ \move{(i)}\
\M{}{\d}{} {\d}{}{} {}{}{\c}
\quad\text{etc.}
$$
{\sc move} (ii)\ One of the $B_i$
moves further into one of the $C_j$,
leaving $A$ unchanged:
$$
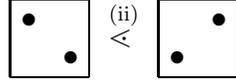
\m{\d}{}{}{\d}\ \move{(ii)}\
\m{}{\d}{\d}{}
$$
{\sc move} (iii)\ The line $A$ lies in $B_i+C_j$.
Then $A$ moves into some $B_{i'}\subset B_i$, 
and so does the corresponding line in $C_j$.
Alternatively, reverse the roles of $B_i$ and
$C_j$.
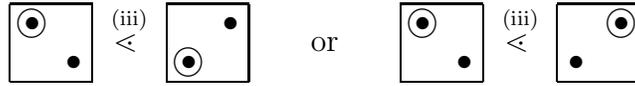
$$
\m{\c}{}{}{\d}\ \move{(iii)}\
\m{}{\d}{\c}{}
\qquad\text{or}\qquad
\m{\c}{}{}{\d}\ \move{(iii)}\
\m{}{\c}{\d}{}
$$
{\sc move} (iv)\ 
The line $A$ lies in 
$B_i+C_j$, but not in $B_{i'}+C_{j'}$,
where $B_{i'}\subset B_i$
and $C_{j'}\subset C_j$.  
Then $A$ moves into $B_{i'}+C_{j'}$,
and the corresponding line in $B_i+C_j$ 
moves with it.
$$
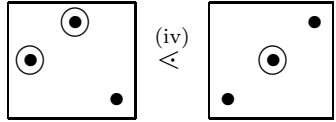
\M{}{\c}{} {\c}{}{} {}{}{\d}
\ \move{(iv)}\
\M{}{}{\d} {}{\c}{} {\d}{}{}
$$
{\sc move} (v)\ The line 
$A$ lies in $B_{i}+C_{j}$. 
Then $B_{i}$ moves further into $C_{j}$,
but $A$ does {\it not} move with it, 
remaining outside $B_{i}\cap C_{j}$.
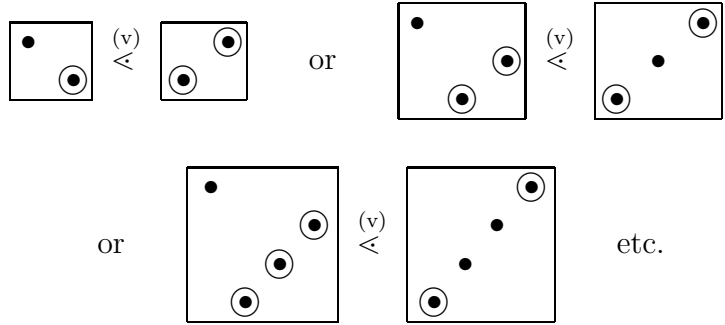
$$
\m{\d}{} {}{\c} \ \move{(v)}\
\m{}{\c} {\c}{}  
\qquad\text{or}\qquad
\M{\d}{}{} {}{}{\c} {}{\c}{}\ \move{(v)}\
\M{}{}{\c} {}{\d}{} {\c}{}{}  
$$
\\[-.9em]
$$
\quad\text{or}\qquad
\begin{array}{|@{\!\ }c@{\!}c@{\!}c@{\!}c@{\!\ }|}
\hline
\d&&& \\ 
&&&\c\\
&&\c&\\
&\c&&\\
\hline
\end{array}
\ \move{(v)}\
\begin{array}{|@{\!\ }c@{\!}c@{\!}c@{\!}c@{\!\ }|}
\hline
&&&\c\\
&&\d&\\
&\d&&\\
\c&&& \\ 
\hline
\end{array}
\qquad\text{etc.}
$$
Note that the underlying 
permutations in this move may differ by an 
arbitrary-length cycle in $S_n$,
not necessarily a transposition.
\\[.5em]
\indent As in the classical case of two flags,
certain regions enclosed by the affected
dots must be empty for these moves to 
define covers $\move{}$
(though they always define
relations $<$). See \S\ref{sec:moves}.
The above moves may seem complicated,
but they are unavoidable in any computationally
effective description: the minimal degenerations
are what they are.

We showed in \cite{mwz1} that the 
number of parameters of
a configuration type (i.e., its dimension
when thought of as a
$\GLn$-orbit in $\PP^{n-1} \sh\times\Flag(\CC^n)\sh\times
\Flag(\CC^n)$ ) is:
$$
\dim(F_{w,\D})=
\tbinom n2+(n\sh-1)+\ell(w)-\#\left\{j\left|\ 
\begin{array}{c}
\text{for all}\
k\in\Delta,\\ k<j\ \text{or}\ w(k)<w(j)
\end{array} \right.\right\}\,.
$$
For example, $F_\min$ has
dimension $\binom n2 + (n\sh-1) + 0 -(n\sh-1)
=\binom n2$.  Indeed, the minimal orbit is isomorphic
to $\Flag(\CC^n)$.  

It is easily seen from the description 
of the moves (i)--(v), 
together with the dot-vanishing
conditions in \S\ref{sec:moves}, that each move
increases the dimension by one.  Thus, 
our Bruhat order is a poset ranked by
$\dim(F)-\dim(F_\min)$.  (This is no longer true
if $(B_\bdot,C_\bdot)$ are partial flags,
and it is not clear whether our poset is ranked.)

We conjecture that 
a refinement of the move-labels
(i)--(v) on the covers of our poset
will give a lexicographic shelling similar to
that of Edelman \cite{edelman}
for (undecorated) permutations.

\subsection{Examples $n=2$, $3$}\label{examples}

We illustrate our constructions in the simplest cases.
Let $n=2$.  Then the Hasse diagram of our Bruhat order is:
\\

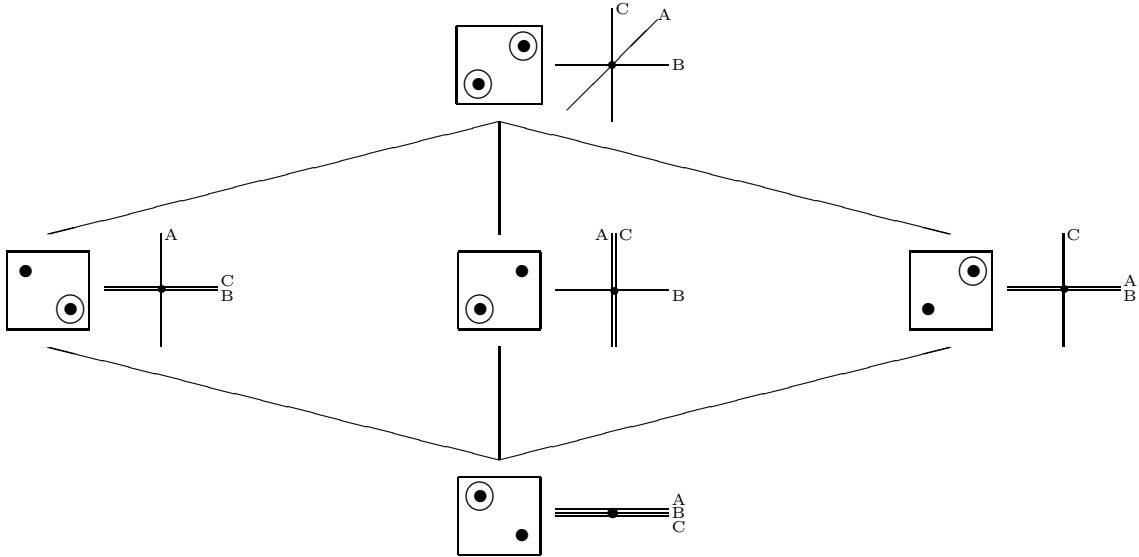
\mbox{}\hspace{-.4in}
\unitlength 1.5mm
\linethickness{0.8pt}
\begin{picture}(60.00,50.00)

\put(40,45){\makebox(0,0)[cc]
{$\m{}{\c}{\c}{}$}}
\put(0,25){\makebox(0,0)[cc]
{$\m{\d}{}{}{\c}$}}
\put(40,25){\makebox(0,0)[cc]
{$\m{}{\d}{\c}{}$}}
\put(80,25){\makebox(0,0)[cc]
{$\m{}{\c}{\d}{}$}}
\put(40,05){\makebox(0,0)[cc]
{$\m{\c}{}{}{\d}$}}

\put(40,10){\line(4,1){40}}
\put(40,10){\line(-4,1){40}}
\put(40,10){\line(0,1){10}}
\put(40,40){\line(4,-1){40}}
\put(40,40){\line(-4,-1){40}}
\put(40,40){\line(0,-1){10}}

\linethickness{0.4pt}

\put(90.07,25.15){\circle*{0.7}}
\put(85,25){\line(1,0){10}}
\put(95,25){\makebox(0,0)[lt]
{\tiny \,B}}
\put(85,25.3){\line(1,0){10}}
\put(95,25.4){\makebox(0,0)[lb]
{\tiny \,A}}
\put(90,20){\line(0,1){10}}
\put(90,30){\makebox(0,0)[lc]
{\tiny \,C}}

\put(50.07,5.3){\circle*{1}}
\put(45,5.0){\line(1,0){10}}
\put(45,5.3){\line(1,0){10}}
\put(45,5.6){\line(1,0){10}}
\put(55,4.5){\makebox(0,0)[lt]
{\tiny \,C}}
\put(55,5.3){\makebox(0,0)[lc]
{\tiny \,B}}
\put(55,6.0){\makebox(0,0)[lb]
{\tiny \,A}}

\put(50.2,25.0){\circle*{0.7}}
\put(45,25){\line(1,0){10}}
\put(55,25){\makebox(0,0)[lt]
{\tiny \,B}}
\put(50.3,20){\line(0,1){10}}
\put(50.3,30){\makebox(0,0)[lc]
{\tiny \,C}}
\put(50,20){\line(0,1){10}}
\put(50,30){\makebox(0,0)[rc]
{\tiny A\,}}

\put(10.1,25.15){\circle*{0.7}}
\put(05,25){\line(1,0){10}}
\put(15,25){\makebox(0,0)[lt]
{\tiny \,B}}
\put(05,25.3){\line(1,0){10}}
\put(15,25.4){\makebox(0,0)[lb]
{\tiny \,C}}
\put(10,20){\line(0,1){10}}
\put(10,30){\makebox(0,0)[lc]
{\tiny \,A}}

\put(50.0,45.0){\circle*{0.7}}
\put(45,45){\line(1,0){10}}
\put(55,45){\makebox(0,0)[lc]
{\tiny \,B}}
\put(50.0,40){\line(0,1){10}}
\put(50.0,50){\makebox(0,0)[lc]
{\tiny \,C}}
\put(46,41){\line(1,1){8}}
\put(54,49){\makebox(0,0)[lb]
{\tiny A\,}}

\end{picture}
\\[1em]
Next to each decorated permutation,
we have sketched the corresponding lines 
$A,\, B\sh=B_1,\,C\sh=C_1$ in $\CC^2$, with {\scriptsize
$\begin{array}{c@{\!}c}
\underline{\underline{\qquad}}&\ A\\[-0.0em]
&\ B\end{array}$}
indicating that $A$ and $B$ coincide.
The elements of our poset correspond to
the $\GL_2$-orbits on $(\PP^1)^3
=(\PP^1)\sh\times\Flag(\CC^2)
\sh\times\Flag(\CC^2)$:
the minimal element is the full diagonal
$\PP^1\subset(\PP^1)^3$; 
the mid-level elements are the
three partial diagonals,
homeomorphic to $\PP^1\sh\times\,\CC$;
and the maximal element is the generic orbit,
homeomorphic to $\PP^1\sh\times
\CC\sh\times\CC^\times$, where $\CC^\times
=\CC\sh\setminus\{\mathop{\rm pt}\}$.
Note that this last orbit is not a topological 
cell, even after fibering out $\PP^1$.

Now let $n=3$.  We can enumerate the configuration
types by counting the possible decorations 
(decreasing subsequences) of each permutation.  
The identity permutation has $n=3$ decorations,
the longest permutation has $2^n-1=7$.
$$
\begin{array}{ccccccccccccc}
3&+&4&+&4&+&5&+&5&+&7&=&28\\
\mbox{\small$\underline{123}$}&&
\mbox{\small$\underline{213}$}&&
\mbox{\small$\underline{132}$}&&
\mbox{\small$\underline{231}$}&&
\mbox{\small$\underline{312}$}&&
\mbox{\small$\underline{321}$}
\end{array}
$$
The Hasse diagram appears on the next page.
We have labelled the elements $\min\!$, $a$, $b$,\ldots,
$x$, $y$, $z$, $\max\!$, as indicated.
For example, 
$$p=(\underline{312},\{1,2\})=
\mbox{\tiny ${\M\n{\sh\c}\n \n\n{\sh\d} {\sh\c}\n\n}$}\ .
$$
The 72 covering relations,
each coming from a move of type (i)--(v),
are:   
$$
\begin{array}{cccccccccc}
\min\!\!\sh\lhd a& \min\!\!\sh\lhd b& \min\!\!\sh\lhd c& \min\!\!\sh\lhd d&
 a \sh\lhd e& a \sh\lhd f& a \sh\lhd g& a \sh\lhd h\\
 b \sh\lhd f& b \sh\lhd g& b \sh\lhd i& b \sh\lhd j& b \sh\lhd k& b \sh\lhd l&
 c \sh\lhd h& c \sh\lhd i& c \sh\lhd l\\
 d \sh\lhd h& d \sh\lhd j& d \sh\lhd k& 
e \sh\lhd m& e \sh\lhd n&
 f \sh\lhd n& f \sh\lhd o& f \sh\lhd q\\
 g \sh\lhd n& g \sh\lhd p& g \sh\lhd r&
 h \sh\lhd m& h \sh\lhd o& h \sh\lhd p& h \sh\lhd q& h \sh\lhd r& h \sh\lhd s\\
 i \sh\lhd p& i \sh\lhd r& i \sh\lhd s& i \sh\lhd u&
 j \sh\lhd o& j \sh\lhd t\\
 k \sh\lhd p& k \sh\lhd q& k \sh\lhd s& k \sh\lhd t&
 l \sh\lhd p& l \sh\lhd u&
 m \sh\lhd v& m \sh\lhd w\\
 n \sh\lhd v& n \sh\lhd w& n \sh\lhd y&
 o \sh\lhd v& o \sh\lhd x& o \sh\lhd y&
 p \sh\lhd w& p \sh\lhd y& p \sh\lhd z\\
 q \sh\lhd w& q \sh\lhd x&
 r \sh\lhd w& r \sh\lhd z&
 s \sh\lhd x& s \sh\lhd z&
 t \sh\lhd x& t \sh\lhd y\\
 u \sh\lhd y& u \sh\lhd z&
 v \sh\lhd \max& w \sh\lhd \max& x \sh\lhd \max& y \sh\lhd \max&\end{array}
$$
The elements in the $i^{\rm th}$ rank of the poset
have orbit dimension 
$i+\dim(F_\min)=i+3$. 



\newcommand{\Min}{\M\c\n\n \n\d\n \n\n\d}
\newcommand{\Ma}{\M\d\n\n \n\c\n \n\n\d}
\newcommand{\Mb}{\M\c\n\n \n\n\d \n\d\n}
\newcommand{\Mc}{\M\n\d\n \c\n\n \n\n\d}
\newcommand{\Md}{\M\n\c\n \d\n\n \n\n\d}
\newcommand{\Me}{\M\d\n\n \n\d\n \n\n\c}
\newcommand{\Mf}{\M\d\n\n \n\n\c \n\d\n}
\newcommand{\Mg}{\M\d\n\n \n\n\d \n\c\n}
\newcommand{\Mh}{\M\n\c\n \c\n\n \n\n\d}
\newcommand{\Mi}{\M\n\n\d \c\n\n \n\d\n}
\newcommand{\Mj}{\M\n\n\c \d\n\n \n\d\n}
\newcommand{\Mk}{\M\n\c\n \n\n\d \d\n\n}
\newcommand{\Ml}{\M\n\d\n \n\n\d \c\n\n}
\newcommand{\Mm}{\M\n\d\n \d\n\n \n\n\c}
\newcommand{\Mn}{\M\d\n\n \n\n\c \n\c\n}
\newcommand{\Mo}{\M\n\n\c \c\n\n \n\d\n}
\newcommand{\Mp}{\M\n\c\n \n\n\d \c\n\n}
\newcommand{\Mq}{\M\n\d\n \n\n\c \d\n\n}
\newcommand{\Mr}{\M\n\n\d \d\n\n \n\c\n}
\newcommand{\Ms}{\M\n\n\d \n\c\n \d\n\n}
\newcommand{\Mt}{\M\n\n\c \n\d\n \d\n\n}
\newcommand{\Mu}{\M\n\n\d \n\d\n \c\n\n}
\newcommand{\Mv}{\M\n\n\c \d\n\n \n\c\n}
\newcommand{\Mw}{\M\n\d\n \n\n\c \c\n\n}
\newcommand{\Mx}{\M\n\n\c \n\c\n \d\n\n}
\newcommand{\My}{\M\n\n\c \n\d\n \c\n\n}
\newcommand{\Mz}{\M\n\n\d \n\c\n \c\n\n}
\newcommand{\Max}{\M\n\n\c \n\c\n \c\n\n}


\newcommand{\Cxmin}{375}
\newcommand{\Cxa}{234}
\newcommand{\Cxb}{328}
\newcommand{\Cxc}{422}
\newcommand{\Cxd}{516}
\newcommand{\Cxe}{47}
\newcommand{\Cxf}{141}
\newcommand{\Cxg}{234}
\newcommand{\Cxh}{328}
\newcommand{\Cxi}{422}
\newcommand{\Cxj}{516}
\newcommand{\Cxk}{609}
\newcommand{\Cxl}{703}
\newcommand{\Cxm}{0}
\newcommand{\Cxn}{94}
\newcommand{\Cxo}{188}
\newcommand{\Cxp}{281}
\newcommand{\Cxq}{375}
\newcommand{\Cxr}{469}
\newcommand{\Cxs}{563}
\newcommand{\Cxt}{656}
\newcommand{\Cxu}{750}
\newcommand{\Cxv}{188}
\newcommand{\Cxw}{281}
\newcommand{\Cxx}{375}
\newcommand{\Cxy}{469}
\newcommand{\Cxz}{563}
\newcommand{\Cxmax}{375}


\newcommand{\Lxmin}{305}
\newcommand{\Lxa}{184}
\newcommand{\Lxb}{278}
\newcommand{\Lxc}{372}
\newcommand{\Lxd}{466}
\newcommand{\Lxe}{-3}
\newcommand{\Lxf}{91}
\newcommand{\Lxg}{184}
\newcommand{\Lxh}{278}
\newcommand{\Lxi}{372}
\newcommand{\Lxj}{466}
\newcommand{\Lxk}{559}
\newcommand{\Lxl}{653}
\newcommand{\Lxm}{-55}
\newcommand{\Lxn}{44}
\newcommand{\Lxo}{138}
\newcommand{\Lxp}{231}
\newcommand{\Lxq}{325}
\newcommand{\Lxr}{419}
\newcommand{\Lxs}{513}
\newcommand{\Lxt}{606}
\newcommand{\Lxu}{700}
\newcommand{\Lxv}{138}
\newcommand{\Lxw}{231}
\newcommand{\Lxx}{325}
\newcommand{\Lxy}{419}
\newcommand{\Lxz}{513}
\newcommand{\Lxmax}{300}


\newcommand{\Cymin}{0}
\newcommand{\Cya}{200}
\newcommand{\Cyb}{200}
\newcommand{\Cyc}{200}
\newcommand{\Cyd}{200}
\newcommand{\Cye}{400}
\newcommand{\Cyf}{400}
\newcommand{\Cyg}{400}
\newcommand{\Cyh}{400}
\newcommand{\Cyi}{400}
\newcommand{\Cyj}{400}
\newcommand{\Cyk}{400}
\newcommand{\Cyl}{400}
\newcommand{\Cym}{600}
\newcommand{\Cyn}{600}
\newcommand{\Cyo}{600}
\newcommand{\Cyp}{600}
\newcommand{\Cyq}{600}
\newcommand{\Cyr}{600}
\newcommand{\Cys}{600}
\newcommand{\Cyt}{600}
\newcommand{\Cyu}{600}
\newcommand{\Cyv}{800}
\newcommand{\Cyw}{800}
\newcommand{\Cyx}{800}
\newcommand{\Cyy}{800}
\newcommand{\Cyz}{800}
\newcommand{\Cymax}{1000}


\newcommand{\Bymin}{-40}
\newcommand{\Bya}{160}
\newcommand{\Byb}{160}
\newcommand{\Byc}{160}
\newcommand{\Byd}{160}
\newcommand{\Bye}{360}
\newcommand{\Byf}{360}
\newcommand{\Byg}{360}
\newcommand{\Byh}{360}
\newcommand{\Byi}{360}
\newcommand{\Byj}{360}
\newcommand{\Byk}{360}
\newcommand{\Byl}{360}
\newcommand{\Bym}{560}
\newcommand{\Byn}{560}
\newcommand{\Byo}{560}
\newcommand{\Byp}{560}
\newcommand{\Byq}{560}
\newcommand{\Byr}{560}
\newcommand{\Bys}{560}
\newcommand{\Byt}{560}
\newcommand{\Byu}{560}
\newcommand{\Byv}{760}
\newcommand{\Byw}{760}
\newcommand{\Byx}{760}
\newcommand{\Byy}{760}
\newcommand{\Byz}{760}
\newcommand{\Bymax}{960}


\pagebreak

\addtolength{\topmargin}{-1.0in}
\addtolength{\oddsidemargin}{-1.0in}
\pagestyle{empty}

\unitlength 0.0095in 
\linethickness{0.8pt}

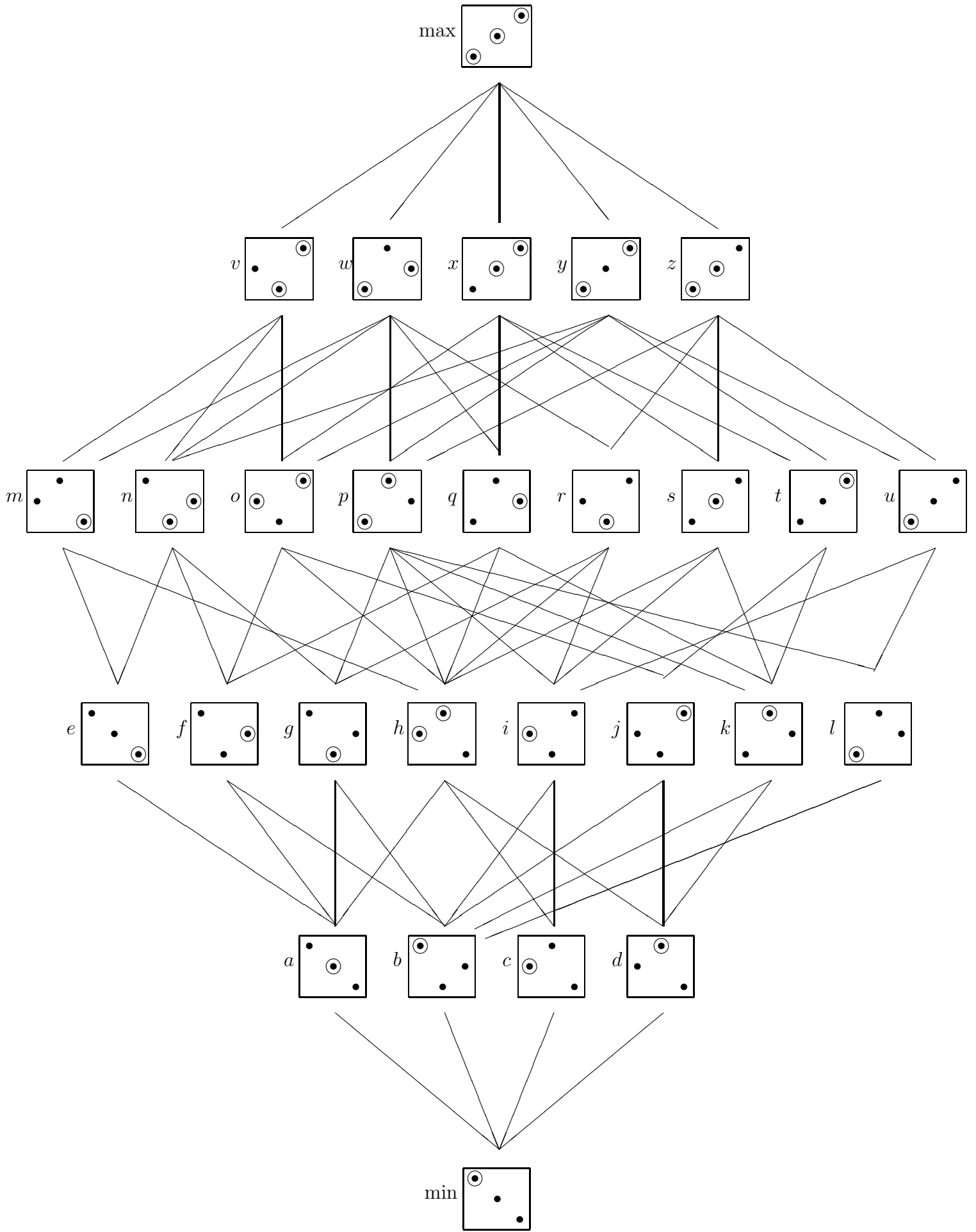
\begin{picture}(750,1000)



\put(\Cxmin,\Cymin){\makebox(0,0)[cc]
 {\footnotesize $\Min$ }}
\put(\Cxa,\Cya){\makebox(0,0)[cc]
 {\footnotesize $\Ma$ }}
\put(\Cxb,\Cyb){\makebox(0,0)[cc]
 {\footnotesize $\Mb$ }}
\put(\Cxc,\Cyc){\makebox(0,0)[cc]
 {\footnotesize $\Mc$ }}
\put(\Cxd,\Cyd){\makebox(0,0)[cc]
 {\footnotesize $\Md$ }}
\put(\Cxe,\Cye){\makebox(0,0)[cc]
 {\footnotesize $\Me$ }}
\put(\Cxf,\Cyf){\makebox(0,0)[cc]
 {\footnotesize $\Mf$ }}
\put(\Cxg,\Cyg){\makebox(0,0)[cc]
 {\footnotesize $\Mg$ }}
\put(\Cxh,\Cyh){\makebox(0,0)[cc]
 {\footnotesize $\Mh$ }}
\put(\Cxi,\Cyi){\makebox(0,0)[cc]
 {\footnotesize $\Mi$ }}
\put(\Cxj,\Cyj){\makebox(0,0)[cc]
 {\footnotesize $\Mj$ }}
\put(\Cxk,\Cyk){\makebox(0,0)[cc]
 {\footnotesize $\Mk$ }}
\put(\Cxl,\Cyl){\makebox(0,0)[cc]
 {\footnotesize $\Ml$ }}
\put(\Cxm,\Cym){\makebox(0,0)[cc]
 {\footnotesize $\Mm$ }}
\put(\Cxn,\Cyn){\makebox(0,0)[cc]
 {\footnotesize $\Mn$ }}
\put(\Cxo,\Cyo){\makebox(0,0)[cc]
 {\footnotesize $\Mo$ }}
\put(\Cxp,\Cyp){\makebox(0,0)[cc]
 {\footnotesize $\Mp$ }}
\put(\Cxq,\Cyq){\makebox(0,0)[cc]
 {\footnotesize $\Mq$ }}
\put(\Cxr,\Cyr){\makebox(0,0)[cc]
 {\footnotesize $\Mr$ }}
\put(\Cxs,\Cys){\makebox(0,0)[cc]
 {\footnotesize $\Ms$ }}
\put(\Cxt,\Cyt){\makebox(0,0)[cc]
 {\footnotesize $\Mt$ }}
\put(\Cxu,\Cyu){\makebox(0,0)[cc]
 {\footnotesize $\Mu$ }}
\put(\Cxv,\Cyv){\makebox(0,0)[cc]
 {\footnotesize $\Mv$ }}
\put(\Cxw,\Cyw){\makebox(0,0)[cc]
 {\footnotesize $\Mw$ }}
\put(\Cxx,\Cyx){\makebox(0,0)[cc]
 {\footnotesize $\Mx$ }}
\put(\Cxy,\Cyy){\makebox(0,0)[cc]
 {\footnotesize $\My$ }}
\put(\Cxz,\Cyz){\makebox(0,0)[cc]
 {\footnotesize $\Mz$ }}
\put(\Cxmax,\Cymax){\makebox(0,0)[cc]
 {\footnotesize $\Max$ }}


\put(\Lxmin,\Cymin){
$\min$}
\put(\Lxa,\Cya){
$a$}
\put(\Lxb,\Cyb){
$b$}
\put(\Lxc,\Cyc){
$c$}
\put(\Lxd,\Cyd){
$d$}
\put(\Lxe,\Cye){
$e$}
\put(\Lxf,\Cyf){
$f$}
\put(\Lxg,\Cyg){
$g$}
\put(\Lxh,\Cyh){
$h$}
\put(\Lxi,\Cyi){
$i$}
\put(\Lxj,\Cyj){
$j$}
\put(\Lxk,\Cyk){
$k$}
\put(\Lxl,\Cyl){
$l$}
\put(\Lxm,\Cym){
$m$}
\put(\Lxn,\Cyn){
$n$}
\put(\Lxo,\Cyo){
$o$}
\put(\Lxp,\Cyp){
$p$}
\put(\Lxq,\Cyq){
$q$}
\put(\Lxr,\Cyr){
$r$}
\put(\Lxs,\Cys){
$s$}
\put(\Lxt,\Cyt){
$t$}
\put(\Lxu,\Cyu){
$u$}
\put(\Lxv,\Cyv){
$v$}
\put(\Lxw,\Cyw){
$w$}
\put(\Lxx,\Cyx){
$x$}
\put(\Lxy,\Cyy){
$y$}
\put(\Lxz,\Cyz){
$z$}
\put(\Lxmax,\Cymax){
$\max$}



\put(\Cxa,\Bya){\line(6, -5){141}} 
\put(\Cxb,\Byb ){\line(2, -5){47}}
\put(\Cxc,\Byc ){\line(-2, -5){47}}
\put(\Cxd,\Byd ){\line(-6, -5){141}}
\put(\Cxe,\Bye ){\line(3, -2){187}}
\put(\Cxf,\Byf ){\line(3, -4){93}}
\put(\Cxg,\Byg ){\line(0, -1){125}}
\put(\Cxh,\Byh ){\line(-3, -4){94}}
\put(\Cxf,\Byf ){\line(3, -2){187}}
\put(\Cxg,\Byg ){\line(3, -4){94}}
\put(\Cxi,\Byi ){\line(-3, -4){94}}
\put(\Cxj,\Byj ){\line(-3, -2){188}}
\put(\Cxk,\Byk ){\line(-2, -1){255}}
\put(\Cxl,\Byl ){\line(-5, -2){340}}
\put(\Cxh,\Byh ){\line(3, -4){94}}
\put(\Cxi,\Byi ){\line(0, -1){125}}
\put(\Cxl,\Byl ){\line(-5, -2){281}}
\put(\Cxh,\Byh ){\line(3, -2){188}}
\put(\Cxj,\Byj ){\line(0, -1){125}}
\put(\Cxk,\Byk ){\line(-3, -4){93}}
\put(\Cxm,\Bym ){\line(2, -5){47}}
\put(\Cxn,\Byn ){\line(-2, -5){47}}
\put(\Cxn,\Byn ){\line(2, -5){47}}
\put(\Cxo,\Byo ){\line(-2, -5){47}}
\put(\Cxq,\Byq ){\line(-2, -1){234}}
\put(\Cxn,\Byn ){\line(6, -5){140}}
\put(\Cxp,\Byp ){\line(-2, -5){47}}
\put(\Cxr,\Byr ){\line(-2, -1){235}}
\put(\Cxm,\Bym ){\line(5, -2){305}}
\put(\Cxo,\Byo ){\line(6, -5){140}}
\put(\Cxp,\Byp ){\line(2, -5){47}}
\put(\Cxq,\Byq ){\line(-2, -5){47}}
\put(\Cxr,\Byr ){\line(-6, -5){141}}
\put(\Cxs,\Bys ){\line(-2, -1){235}}
\put(\Cxp,\Byp ){\line(6, -5){141}}
\put(\Cxr,\Byr ){\line(-2, -5){47}}
\put(\Cxs,\Bys ){\line(-6, -5){141}}
\put(\Cxu,\Byu ){\line(-5, -2){305}}
\put(\Cxo,\Byo ){\line(3, -1){328}}
\put(\Cxt,\Byt ){\line(-5, -4){140}}
\put(\Cxp,\Byp ){\line(5, -2){305}}
\put(\Cxq,\Byq ){\line(2, -1){234}}
\put(\Cxs,\Bys ){\line(2, -5){46}}
\put(\Cxt,\Byt ){\line(-2, -5){47}}
\put(\Cxp,\Byp ){\line(4, -1){417}}
\put(\Cxu,\Byu ){\line(-1, -2){53}}
\put(\Cxv,\Byv ){\line(-3, -2){188}}
\put(\Cxw,\Byw ){\line(-2, -1){250}}
\put(\Cxv,\Byv ){\line(-4, -5){100}}
\put(\Cxw,\Byw ){\line(-3, -2){187}}
\put(\Cxy,\Byy ){\line(-3, -1){375}}
\put(\Cxv,\Byv ){\line(0, -1){125}}
\put(\Cxx,\Byx ){\line(-3, -2){187}}
\put(\Cxy,\Byy ){\line(-2, -1){250}}
\put(\Cxw,\Byw ){\line(0, -1){125}}
\put(\Cxy,\Byy ){\line(-3, -2){188}}
\put(\Cxz,\Byz ){\line(-2, -1){250}}
\put(\Cxw,\Byw ){\line(4, -5){94}}
\put(\Cxx,\Byx ){\line(0, -1){120}}
\put(\Cxw,\Byw ){\line(5, -3){188}}
\put(\Cxz,\Byz ){\line(-4, -5){92}}
\put(\Cxx,\Byx ){\line(3, -2){188}}
\put(\Cxz,\Byz ){\line(0, -1){125}}
\put(\Cxx,\Byx ){\line(2, -1){250}}
\put(\Cxy,\Byy ){\line(3, -2){187}}
\put(\Cxy,\Byy ){\line(2, -1){250}}
\put(\Cxz,\Byz ){\line(3, -2){187}}
\put(\Cxmax,\Bymax ){\line(-3, -2){187}}
\put(\Cxmax,\Bymax ){\line(-4, -5){94}}
\put(\Cxmax,\Bymax ){\line(0, -1){120}}
\put(\Cxmax,\Bymax ){\line(4, -5){94}}
\put(\Cxmax,\Bymax ){\line(3, -2){188}}

\end{picture}

\pagebreak

\addtolength{\topmargin}{1.0in}
\addtolength{\oddsidemargin}{1.0in}
\pagestyle{plain}

As an illustration of the tableau criterion 
(i.e., rank numbers defining the Bruhat order), 
let us check that $e\not\leq z$: that is,
$$
\mbox{
{\tiny $\begin{array}{|@{\!}c@{\!}c@{\!}c@{\!}|}\hline 
 \ \bullet\ &&  \\
 &\ \bullet\ &\\
& & \cbullet\\
\hline\end{array}
$}
$=(\underline{123},\{3\})
=(A,B_\bdot,C_\bdot)$
\ \ and\ \ \
{\tiny $\begin{array}{|@{\!}c@{\!}c@{\!}c@{\!}|}\hline 
& &\ \bullet\ \\
 &\cbullet&\\
  \cbullet &&  \\
\hline\end{array}
$}
$=(\underline{321},\{1,2\})=(A',B'_\bdot,C'_\bdot)$
}
$$
are unrelated elements in our poset,
even though $\underline{123}<\underline{321}$
in the classical Bruhat order.
Indeed, in the second configuration,
$A'=\langle v_1+v_2\rangle \subset
B'_2=\langle v_1,v_2\rangle$,
and no degeneration of $(A',B'_\bdot,C'_\bdot)$
can destroy this containment.
However, in the first configuration,
$A=\langle v_3\rangle \not\subset
B_2=\langle v_1,v_2\rangle$.  Thus
$(\underline{123},\{3\})
\not\leq(\underline{321},\{1,2\})$.
In terms of our rank numbers:
$r_{ij}(\underline{123})\geq
r_{ij}(\underline{321})$ for all $i,j$,
in particular $r_{11}(\underline{123})>
r_{11}(\underline{321})$; but
$\rb{20}(\underline{123},\{3\})<
\rb{20}(\underline{321},\{1,2\})$.

\subsection{Structure of the paper}\label{structure}

Now we sketch our proof of the above results.
After some easy geometric arguments,
we reduce our claims to a rather difficult
combinatorial lemma.
The idea is to approximate the geometric
degeneration order from
above and below by combinatorially defined orders,
and then show that these combinatorial
bounds are equal.

To begin, we distinguish in \S\ref{three}
three partial orders
on decorated permutations $(w,\Delta)$.
First, our geometric order
$\ldeg$ defined by degenerations
of the corresponding configuration
types $F_{w,\Delta}$.
Second, the combinatorial order
$\lrk$ defined in terms of the rank
numbers $r_{ij}(w)$,\, $\rb{ij}(w,\Delta)$.
Third, the order $\lmv$ generated by repeated 
application of our moves $\move{i}$,\ldots,
$\move{v}$.
We wish to show the equivalence of these
three orders.

Some simple geometry and linear algebra
in \S\ref{geom} suffices to show that:
$$
(w,\Delta)\lmv(y,\Gamma)
\,\Rightarrow\,
(w,\Delta)\ldeg(y,\Gamma)
\,\Rightarrow\,
(w,\Delta)\lrk(y,\Gamma)\,.
$$
That is, any move corresponds to a degeneration, and any degeneration increases the rank numbers.
We are then left in \S\ref{combin} to show 
the purely combinatorial
assertion: $$(w,\Delta)\lrk(y,\Gamma)
\,\Rightarrow\,
(w,\Delta)\lmv(y,\Gamma).$$
Given a relation
$(w,\Delta)\slrk(y,\Gamma)$,
we find a move $(w,\Delta)\move{mv} (\tilde w,\tilde\Delta)$ such that the smaller rank numbers of 
$(\tilde w,\tilde\Delta)$ still dominate those
of $(y,\Gamma)$:
$$
(w,\Delta)\move{mv}(\tilde w,\tilde\Delta)
\lrk (y,\Gamma)\,.
$$
Iterating this construction within our finite
poset, we eventually get 
$$
(w,\Delta)
\move{mv}(\tilde w_1,\tilde\Delta_1)
\move{mv}\cdots
\move{mv}(\tilde w_k,\tilde\Delta_k)
= (y,\Gamma)\,.
$$

Throughout our proof, 
we work in the more general case
where $B_\bdot, C_\bdot$ are arbitrary
partial flags, with orbits indexed not by permutations
but by double cosets of permutations, or ``transport
matrices'', as defined in \S\ref{sec:abruhat1}.
Also, our proofs are characteristic-free: 
our vector spaces are over an 
arbitrary infinite field $\kk$, 
not necessarily $\CC$.

Our Rank Theorem, giving the equivalence of
$\ldeg$ and $\lrk$, is a strengthened converse to
Prop.~4.5 in our work \cite{mwz1}, which relied
heavily on quiver theory.
In the notation of \cite{mwz1},
for a triple flag $X=(A,B_\bdot,C_\bdot)$,
we have:  
$r_{ij}(X)=
\dim\mathop{\rm Hom}(I_{\{(i,j)\}}, X)$
and 
$\rb{ij}(X)=
\dim\mathop{\rm Hom}(I_{\{(i,r),(q,j)\}}, X)$.

Our approach is closely related to that of 
Zwara and Skowronski \cite{zwara1}, \cite{zwara2}; 
Bongartz \cite{bongartz}, Riedtmann \cite{riedtmann},
et~al., who considered the degeneration order on
quiver representations.
However, in our special case, 
our results are sharper than those of the 
general theory.
Our description of the covering moves (i)-(v) 
can be deduced (with non-trivial work)
from Zwara's results on the extension order $\smallstack{ext}{\leq}$ in \cite{zwara1}, \cite{zwara2}.  
But our results about the
rank order $\lrk$ are considerably stronger 
than any of the corresponding general results: our
order requires computing $\mathop{\rm Hom}$ with only a few
indecomposables, rather than all.
\\[1em]
{\bf Acknowlegement}\ \  This work grew out of
a project with Jerzy Weyman and Andrei Zelevinsky.  The author is indebted to
them for suggesting the topic, as well
as numerous helpful suggestions.

\section{Results}\label{three}

\subsection{Two Flags}\label{sec:abruhat1}

In order to establish the notation for our
main theorem in its full generality, 
we first state the classical
theory for two flags.

Throughout this paper, all vector spaces are
over a fixed field $\kk$ of arbitrary 
characteristic, infinite but not necessarily algebraically closed; and we fix
a vector space $V$ of dimension $n$ with 
standard basis $e_1,\ldots,e_n$.  
Let $\bb = (b_1, \ldots, b_q)$ 
be a  list of positive integers 
with sum equal to $n$:  
that is, a \emph{composition} of $n$.
We denote by $\Flag(\bb)$ the variety of partial flags
$B_\bdot = (0 \sh= B_0 \sh\subset B_1 \sh\subset \cdots 
\sh\subset B_q \sh= V)$
of vector subspaces in $V$ such that
$$\dim (B_{i}/B_{i-1}) = b_i \quad (i = 1, \ldots, q) \ .$$
$\Flag(\bb)$ is a homogeneous space under the
natural action of the general linear group 
$\GL(V)=\GLn(\kk)$.

Let us fix two compositions of $n$,\ \,$\bb \sh= (b_1, \ldots, b_q)$
and  $\cc\sh=(c_1,\ldots,c_r)$.  The Schubert (or Bruhat)
decomposition classifies the orbits of $GL(V)$ acting diagonally on 
the double flag variety $\Flag(\bb)\times\Flag(\cc)$.
We index these orbits by {\it transport matrices}
$M=(m_{ij})$, which are $q\sh\times r$ matrices
of nonnegative integers $m_{ij}$ with row sums $b_i=\sum_j m_{ij}$
and column sums $c_j=\sum_i m_{ij}$ (so that the sum of all
entries is $n$). 
If $\bb=\cc=(1,\ldots,1)=(1^n)$, 
then $\Flag(\bb)\sh\times\Flag(\cc)$ consists of pairs 
of full flags, and each transport matrix
is the permutation matrix $M=M_w$
corresponding to a $w\in S_n$,
with $m_{w(i),i}\sh=1$ and $m_{ij}\sh=0$ otherwise.

Given a transport matrix $M$,
we define the orbit $F_M \subset\Flag(\bb)\sh\times\Flag(\cc)$
as the following set of double flags $(B_\bdot,C_\bdot)$.
Given any basis of $V$ with the $n$ vectors indexed as:
$$
V=\left\langle v_{ijk}\left|\begin{array}{c}
(i,j)\in[q]\sh\times[r]\\ 
1\sh\leq k\sh\leq m_{ij}
\end{array}\right.\!\!\!\right\rangle,
$$
where $[q]=[1,q]:=\{1,2,\ldots,q\}$, let
$$
B_i:=\langle v_{i'jk}\mid i'\leq i\rangle,
\qquad
C_j:=\langle v_{ij'k}\mid j'\leq j\rangle\,,
$$
where $\langle\ \rangle$ denotes linear span.
As the basis $\langle v_{ijk}\rangle$ varies, $(B_\bdot,C_\bdot)$ runs over
all double flags in $F_M$. 
In the case $\bb=\cc=(1^n)$,
with $M=M_w$, we may take
$v_{ij1}=v_i$ for any basis $v_1,\ldots,v_n$
of $V$, and obtain the configuration type $F_w$
from the Introduction. 

We can also describe this orbit by 
intersection conditions:
$$
F_M=\{\,(B_\bdot,C_\bdot)\mid \dim(B_i\cap C_j)=r_{ij}(M)\,\}\,,
$$
where 
$$
r_{ij}(M):=
\subsub{\sum}{(k,l)}{k\leq i,\,l\leq j}\!\!\! m_{kl}
$$
are the {\it rank numbers}.  
This characterization follows from Theorem
\ref{Abruhat} below.

These orbits cover the double flag variety:
$$
\Flag(\bb)\sh\times\Flag(\cc)
=\coprod_M F_M
$$
where the union runs over all transport matrices $M$.

We shall need the following partial order on the matrix positions 
$(i,j)\in[1,q]\sh\times[1,r]$:  we write 
$$
(i,j)\leq(i',j')  \quad\Longleftrightarrow\quad
i\sh\leq i'\ \text{and}\ j\sh\leq j'\,.
$$
That is, the northwest positions are small, the 
southeast positions large.
Also, $(i,j)<(i',j')$ means $(i,j)\leq(i',j')$ and $(i,j)\neq(i',j')$,
a convention we will use when dealing with any partial order.
Furthermore, for sets of positions $\Delta,\Delta'\subset[1,q]\sh\times[1,r]$
we let:
$$
\Delta\leq\Delta'  \quad\Longleftrightarrow\quad
\forall\,(i,j)\sh\in\Delta\ \ \exists\,(i',j')\sh\in\Delta'\
\text{with}\ (i,j)\leq(i',j')\,.
$$

Now, the {\it degeneration order} or {\it Ehresmann-Bruhat order}
on the set of all transport matrices
describes how the orbits $F_M$ touch each other:
$$
M\ldeg M' \quad\Longleftrightarrow\quad  
\overline{F}_M\subset \overline{F}_{M'}
\quad\Longleftrightarrow\quad  
F_M\subset\overline{F}_{M'}\ ,
$$
where $\overline F_M$ denotes the (Zariski) closure of $F_M$.
Our goal is to give a combinatorial characterization of this geometric
order.  

First, we approximate the degeneration order on double flags 
by comparing rank numbers. We define: 
$$
M\lrk M' \quad\Longleftrightarrow\quad  
r_{ij}(M)\geq r_{ij}(M')\quad
\forall\,(i,j)\in[1,q]\sh\times[1,r]\,.
$$

Second, we define certain moves on matrices 
which will turn out to be the covers
of the degeneration order:  that is, the 
relations $M\move{deg} M'$ such that $M\sldeg M''\ldeg M'\,\Rightarrow\, M''\sh=M'$.
Suppose we consider positions $(i_0,j_0)\sh\leq(i_1,j_1)$
defining a rectangle 
$$
R=[i_0,i_1]\sh\times[j_0,j_1]\subset [1,q]\sh\times[1,r]\,,
$$
and we are given an $M$ satisfying $m_{i_0j_0},m_{i_1j_1}>0$
and $m_{ij}=0$ for all $(i,j)\in R$,
$(i,j)\neq(i_0,j_0),(i_1,j_1),(i_0,j_1),(i_1,j_0)$.  Then
the {\it simple move} on the matrix $M$, at the rectangle $R$,
is the operation which produces the new matrix:
$$
M'=M-E_{i_0j_j}-E_{i_1j_1}+E_{i_0j_1}+E_{i_1j_0}\,,
$$
where $E_{ij}$ denotes the coordinate matrix with 1 in position $(i,j)$
and 0 elsewhere.  We write $M\move{R} M'$ or 
$M\move{mv} M'$. Pictorially,
$$
M=\hspace{-1em}\begin{array}{c@{\!}l}
&\hspace{2em}j_0\hspace{5em}j_1\\
\begin{array}{c}   i_0\\ \\ \\ \\ i_1\\ \end{array}
&\begin{array}{|c@{\!}ccccc@{\!}c|}
\hline
&\vdots&&&&\vdots&\\
\cdots\ & a& 0&\cdots & 0&b &\ \cdots\\
&0&&&&0&\\
&\vdots&&0&&\vdots&\\
&0&&&&0&\\
\cdots\ &c&0&\cdots&0&d&\ \cdots\\
&\vdots&&&&\vdots&\\
\hline
\end{array}
\end{array}
\ \move{R}\  
M'=\hspace{-1em}\begin{array}{c@{\!}l}
&\hspace{2.7em}j_0\hspace{6.3em}j_1\\
\begin{array}{c}   i_0\\ \\ \\ \\ i_1\\ \end{array}
&\begin{array}{|c@{\!}ccccc@{\!}c|}
\hline
&\vdots&&&&\vdots&\\
\cdots\ & a\sh-1 & 0&\cdots & 0& b\sh+1 &\ \cdots\\
&0&&&&0&\\
&\vdots&&0&&\vdots&\\
&0&&&&0&\\
\cdots\ & c\sh+1 &0&\cdots&0& d\sh-1 &\ \cdots\\
&\vdots&&&&\vdots&\\
\hline
\end{array}
\end{array}.
$$

In the case where $M$ is a permutation matrix, 
the simple move corresponds to multiplying by a transposition:
if $M$ is associated to $w$, then $M'$ is associated to 
$w'=(i_0,i_1)\cdot w=w\cdot(j_0,j_1)$, and the vanishing conditions
on $m_{ij}$ assure that the Bruhat length $\ell(w')=\ell(w)+1$.

We say $M$, $M''$ are related by the {\it move order},
if, starting with $M$, we can perform a sequence
of simple moves on various rectangles $R_1, R_2,\ldots$ to obtain $M''$:
$$
M\lmv M''\quad\Longleftrightarrow\quad
M\move{R_1}M'\move{R_2}\cdots M''.
$$

\begin{thm}{(Ehresmann-Chevalley)}\label{Abruhat}\\
(a) The three orders defined above are equivalent:
$$\quad M\ldeg M'\ \ \Longleftrightarrow\ \
M\lrk M'\ \ \Longleftrightarrow\ \
M\lmv M'.$$
(b) The relation $M\sldeg M'$ is a cover exactly
when $M\move{\rm mv} M'$.
\end{thm}

We give a proof in \S\ref{Abruhat proof}.
Once we have established the equivalences above, 
we call the common order the {\it Bruhat order},
written simply as $M\leq M'$.
 
\subsection{A Line and Two Flags}\label{decorated}

We now state our main theorems.
We consider $\GL(V)$ acting diagonally on
$
\PP(V)\times\Flag(\bb)\times\Flag(\cc),
$
the variety of triples of a line and two flags.
We showed in \cite{mwz1} 
that the orbits correspond
to the {\it decorated matrices} $(M,\Delta)$, 
meaning that $M$ is a transport matrix, and 
$$
\Delta=\{(i_1,j_1),\ldots,(i_t,j_t)\}\subset[1,q]\sh\times[1,r]
$$ 
is a set of matrix positions satisfying:
$$
i_1<i_2<\cdots<i_t,\ \ j_1>j_2>\cdots>j_t,\quad\text{and}\quad
m_{ij}>0\ \,\forall\,(i,j)\in\Delta.
$$
That is, the positions  $(i_1,j_1),(i_2,j_2),\ldots$
proceed from northeast to southwest.  
We may concisely write down $(M,\Delta)$ by drawing a circle around
the nonzero entries of $M$ at the positions $(i,j)\in\Delta$.

The corresponding orbit
$F_{M,\Delta}$ consists of the triple flags
$(A,B_\bdot,C_\bdot)$ defined as follows.
Given a basis $\langle v_{ijk}\rangle$ as above,
the flags $B_\bdot, C_\bdot$ are defined exactly as before
(and thus depend only on $M$); and the line is defined as
$A:=\langle\, \sum_{(i,j)\in\Delta} v_{ij1}\,\rangle$.
Thus $M$ indicates the relative positions 
of the two flags $B_\bdot,C_\bdot$,
and $\Delta$ is a ``decoration'' on $M$ 
indicating the position of the line $A$.
Once again we have:
$$
\PP(V)\sh\times\Flag(\bb)\sh\times\Flag(\cc)
=\coprod_{(M,\Delta)} F_{M,\Delta},
$$
where $(M,\Delta)$ runs over all decorated matrices.
We also define the degeneration order $(M,\Delta)\ldeg(M',\Delta')$
as before. 

Next we define the rank order.
For $(A,B_\bdot,C_\bdot)\in F_{M,\Delta}$,
we define a new rank number:
$$
\begin{array}{rcl}
\rb{ij}(M,\Delta)&:=& \dim(B_i\cap C_j)+\dim(A\cap(B_i\sh+C_j))\\[.5em]
&=&r_{ij}(M)+\delta_{ij}(\Delta)\,;
\end{array}
$$
where we define:
$$\begin{array}{rcl}
\delta_{ij}(\Delta)
&:=&\left\{\begin{array}{cl}
1 & \text{if}\ \Delta\leq\{(i,r),(q,j)\}\\
0 & \text{otherwise}\\
\end{array}\right.\\[1em]
&=&\left\{\begin{array}{cl}
0 & \text{if}\ (i\sh+1,j\sh+1)\leq\Delta\\
1 & \text{otherwise.}
\end{array}\right.
\end{array}$$
We can extend this definition to $(i,j)\in[0,q]\sh\times[0,r]$
by setting $r_{ij}(M):=0$ if $i\sh=0$ or $j\sh=0$,
so that $$
\rb{i0}(M,\D)\sh=\dim(A\cap B_i)\sh=
\delta_{i0}(\D),
\quad
\rb{0j}(M,\D)\sh=\dim(A\cap C_j)\sh=
\delta_{0j}(\D).
$$
Now we let:
$$
(M,\Delta)\lrk(M',\Delta')\ \ \Longleftrightarrow\!
\begin{array}{c}
r_{ij}(M)\geq r_{ij}(M')\\
\rb{ij}(M,\D)\geq\rb{ij}(M',\D')
\end{array}\ \forall\,(i,j)\in[0,q]\sh\times[0,r].
$$
\\[.5em]   
{\bf Rank Theorem}\ \ {\it
\label{Sbruhat}
The degeneration order and the rank order are equivalent:
$$
(M,\D)\ldeg (M',\D')\quad\Longleftrightarrow\quad
(M,\D)\lrk(M',\D').
$$
That is, the triple flag $(A,B_\bdot,C_\bdot)$
is a degeneration of $(A',B'_\bdot,C'_\bdot)$
if and only if, for all $(i,j)\in[0,q]\sh\times[0,r]$,
$$
\begin{array}{c}
\dim(B_i\cap C_j)\geq\dim(B'_i\cap C'_j)
\\[.3em]
\dim(B_i\cap C_j)+\dim(A\cap(B_i\sh+C_j))\geq
\dim(B'_i\cap C'_j)+\dim(A'\cap(B'_i\sh+C'_j))\,.
\end{array}
$$
}
   
\subsection{Simple Moves}\label{sec:moves}
\label{moves}

Below, we define simple moves of types (i)--(v)
on a decorated matrix $(M,\D)$, each producing
a new matrix $(M',\D')$, so that we write
$(M,\D)\move{mv}(M',\D')$.  Given these moves, we
define the move order $(M,\D)\lmv(M'',\D'')$ as before.
\\[1em]
{\bf Move Theorem}\ \ {\it
The degeneracy order and the move order are equivalent:
$$
(M,\D)\ldeg(M',\D')\quad\Longleftrightarrow\quad(M,\D)\lmv(M',\D')\,.
$$
}
Again, we call the common order the {\it Bruhat order}.
\\[1em]
{\bf Minimality Theorem}\ \ {\it
The relation $(M,\D)\sldeg(M',\D')$ is a 
cover exactly when $(M,\D)\move{\rm mv}(M',\D')$
for one of the simple moves (i)--(v).
}\\[1em]
We introduce an operation which normalizes an arbitrary subset $S$
of matrix positions into a decoration $\Delta$ of the prescribed form.
For $S\subset [1,q]\sh\times[1,r]$,
let: 
$$
[S]:=\{\,(i,j)\in S\mid  (i,j)\not< (k,l)\ \forall\,(k,l)\sh\in S\,\}
$$
be the set of $\leq$-maximal positions in $S$.
This operation is ``explained'' by the following lemma, proved in
\S\ref{combin}:

\begin{lemma}[Uncircling lemma]
\label{uncirc}
If $M$ is a transport matrix, $S$ a set of matrix positions
with $m_{ij}>0$ for all $(i,j)\in S$, and we define 
$(A,B_\bdot,C_\bdot)$ by the same formulas as for a decorated
matrix (namely $A:=\langle\, \sum_{(i,j)\in S} v_{ij1}\,\rangle$\,),
then $(A,B_\bdot,C_\bdot)\in F_{M,\Delta}$, where $\Delta=[S]$.
\end{lemma}

It remains to define the five types of simple moves
$(M,\D)\!\move{(i)}\!(M',\D')$,\ldots, $(M,\D)\!\move{(v)}\!(M',\D')$.
Although geometrically, it is natural to 
think of the more general configuration
degenerating to the more special
one, combinatorially
it is more convenient to describe the modification
of the more special (smaller) element $(M,\D)$ 
to obtain the more general (larger) element $(M',\D')$.  

In each case, 
we indicate the matrix positions in $(M,\D)$
modified by the move, 
and we list the requirements 
on $M=(m_{ij})$ and $\D$ for the move
to be valid.  
Then we specify the result of the move, 
$(M',\D')$.
\\[.3em]
(i)\ Suppose we have a position $(i_1,j_1)\not\leq\D$,
with $m_{i_1j_1}>0$ and $m_{ij}=0$ whenever 
$(i,j)\sh<(i_1,j_1)$,\, $(i,j)\!\not\leq\!\D$.
Then define:
$$
M':=M,\qquad \D':=[\,\D\cup\{(i_1,j_1)\}\,]\,.
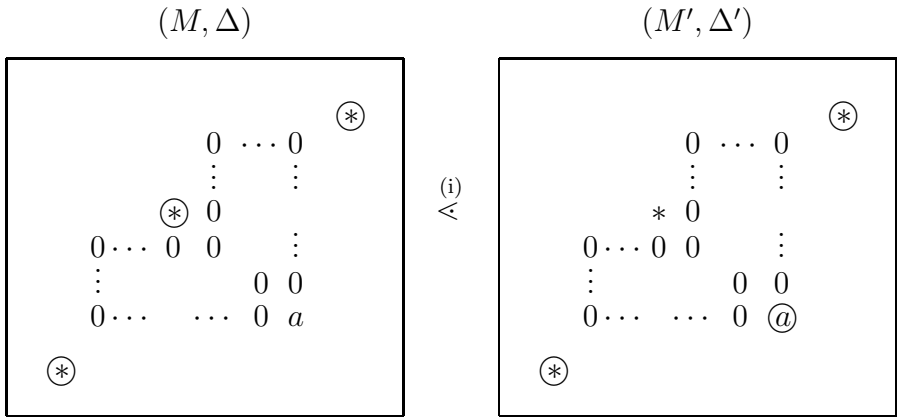
$$
$$\hspace{-.8in}
\begin{array}{c}(M,\D)\\[.5em]
\begin{array}{|c@{\!\ }c@{\!\ }c@{\!\ }
c@{\!\ }c@{\!\ }c@{\!\ }c@{\!\ }c@{\!\ }c@{\!\ }c@{\!\ }c@{\!\ }c@{\!\ }c|}\hline
\ \ &&&&&&&&&&&&\ \ \\
&&&&&&&&&&\ \ &\C{\ast}&\\[-.3em]
&&&&&&&0&\cdots&0&&&\\[-.5em]
&&&&&&&\vdots&&\vdots&&&\\[-.2em]
&&&&&&\C{\ast}&0&&&&&\\[-.5em]
&&&&0&\cdots&0&0&&\vdots&&&\\[-.5em]
&&&&\vdots&&&&0&0&&&\\[-.2em]
&&&&0&\cdots&&\cdots&0&a&&&\\[.5em]
&\C{\ast}&&&&&&&&&&&\\[.8em]
\hline\end{array}
\end{array}
\ \move{(i)}\
\begin{array}{c}(M',\D')\\[.5em]
\begin{array}{|c@{\!\ }c@{\!\ }c@{\!\ }
c@{\!\ }c@{\!\ }c@{\!\ }c@{\!\ }c@{\!\ }c@{\!\ }c@{\!\ }c@{\!\ }c@{\!\ }c|}\hline
\ \ &&&&&&&&&&&&\ \ \\
&&&&&&&&&&\ \ &\C{\ast}&\\[-.3em]
&&&&&&&0&\cdots&0&&&\\[-.5em]
&&&&&&&\vdots&&\vdots&&&\\[-.2em]
&&&&&&{\ast}&0&&&&&\\[-.5em]
&&&&0&\cdots&0&0&&\vdots&&&\\[-.5em]
&&&&\vdots&&&&0&0&&&\\[-.2em]
&&&&0&\cdots&&\cdots&0&\C{a}&&&\\[.5em]
&\C{\ast}&&&&&&&&&&&\\[.8em]
\hline\end{array}
\end{array}
$$
\\
where $a\sh=m_{i_1j_1}\sh>0$,\ the symbol $\ast$ represents a matrix entry
$m_{ij}\sh>0$, and a circle $\C{\ast}$ around $\ast\sh=m_{ij}$ means that 
$(i,j)\sh\in\D$.  The values of $m_{ij}$ are not changed by the move.
\\[1em]
(ii) Suppose we have positions $(i_0,j_0)<(i_1,j_1)$ with:
$m_{i_0j_0},m_{i_1j_1}\sh>0$, and $m_{ij}=0$ whenever:
$(i_0,j_0)\sh<(i,j)\sh<(i_1,j_1)$,\, $(i,j)\neq(i_0,j_1),(i_1,j_0)$.
Suppose further that  $(i_1,j_1)\!\not\in\!\D$;
\ $(i_0,j_1)\!\not\in\!\D$ or $(i_1,j_0)\!\not\in\!\D$; and
$(i_0,j_0)\not\in\Delta$ or $m_{i_0j_0}>1$.
Then define:
$$
M':=M-E_{i_0j_0}-E_{i_1j_1}+E_{i_0j_1}+E_{i_1j_0},\qquad \D':=\D\,.
$$
$$
\begin{array}{c} (M,\D)\\[.5em]
\begin{array}{|c@{\!}ccccc@{\!}c|}
\hline
\ \ \ &&&&&&\ \ \ \\
& a & 0&\cdots & 0& b&\\
&0&&&&0&\\[-.5em]
&\vdots&&0&&\vdots&\\
&0&&&&0&\\
& c&0&\cdots&0& d&\\[1em]
\hline
\end{array}
\end{array}
\ \move{(ii)}\  
\begin{array}{c} (M',\D')\\[.5em]
\begin{array}{|c@{\!}ccccc@{\!}c|}
\hline
\ \ \ &&&&&&\ \ \ \\
& a\!\!-\!\!1 & 0&\cdots & 0& b\!\!+\!\!1 &\\
&0&&&&0&\\[-.5em]
&\vdots&&0&&\vdots&\\
&0&&&&0&\\
& c\!\!+\!\!1 &0&\cdots&0& d\!\!-\!\!1 &\\[1em]
\hline
\end{array}
\end{array}
$$
\\
Here $a\sh=m_{i_0j_0}\sh>0$,\, $d\sh=m_{i_1j_1}\sh>0$,\
$a$ may be circled only if $a\sh>1$, at most one of
$b,c$ is circled, and $d$ is not circled.
The positions of circles are unchanged by the move.
\\[1em]
(iii)(a)  Suppose $(i_0,j_0)\sh<(i_1,j_1)$
with  $(i_0,j_0)\sh\in\D$,\, $m_{i_0j_0}\sh=1$,\ $m_{i_1j_1}\sh>0$,
and $m_{ij}\sh=0$ whenever:
$(i_0,j_0)\sh<(i,j)\sh<(i_1,j_1)$,\, $(i,j)\sh\neq(i_1,j_0)$; 
and whenever $(i,j)\sh\leq(i_0,j_1)$, $(i,j)\!\not\leq\!\D$.
Then define:
$$
M':=M-E_{i_0j_0}-E_{i_1j_1}+E_{i_0j_1}+E_{i_1j_0},\qquad 
\D':=[\,\D\cup\{(i_0,j_1)\}\,]\,.
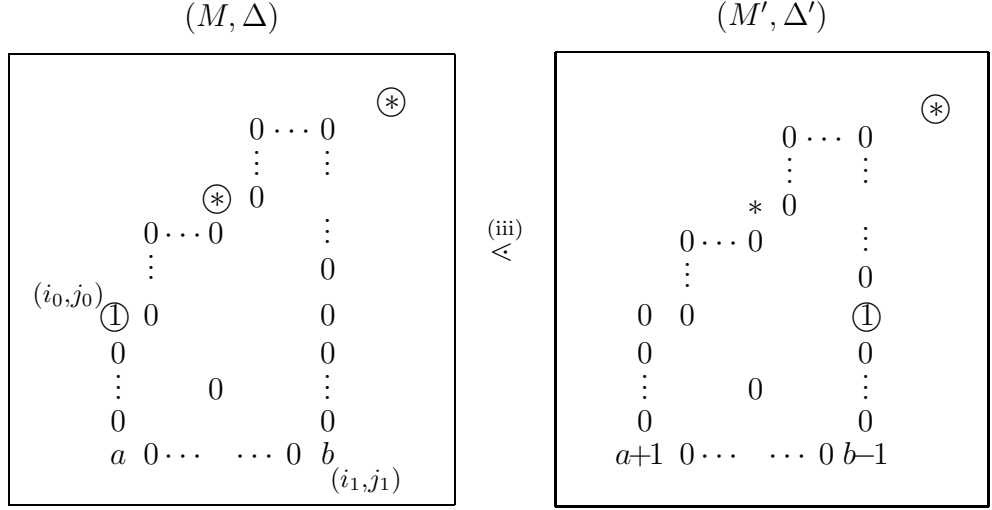
$$
$$\hspace{-0in}
\begin{array}{c}(M,\D)\\[.5em]
\begin{array}{|c@{\!\ }c@{\!\ }c@{\!\ }c@{\!}c@{\!\ }c@{\!}
c@{\!\ }c@{\!\ }c@{\!\ }c@{\!\ }c|}\hline
\ \ &&&&&&&&&&\ \ \ \\[-.3em]
&&&&&&&&\ \ &\C{\ast}&\\[-.3em]
&&&&&0&\cdots&0&&&\\[-.5em]
&&&&&\vdots&&\vdots&&&\\[-.2em]
&&&&\C{\ast}&0&&&&&\\[-.5em]
&&0&\cdots&0&&&\vdots&&&\\[-.5em]
&&\vdots&&&&&0&&&\\[-.5em]
\ \mbox{\footnotesize $(i_0,\!j_0)$}\!\!&&&&&&&&&&\\[-.5em]
&\C{1}\ &0&&&&&0&&&\\[-0em]
&0&&&&&&0&&&\\[-.5em]
&\vdots&&&0&&&\vdots&&&\\[-.2em]
&0&&&&&&0&&&\\[-.2em]
&a&0&\cdots&&\cdots&0&b&&&\\[-.5em]
&&&&&&&&\!\!\mbox{\footnotesize $(i_1,\!j_1)$}\hspace{-1em}&&\\[.2em]
\hline\end{array}
\end{array}
\ \move{(iii)}\
\begin{array}{c}(M',\D')\\[.5em]
\begin{array}{|c@{\!\ }c@{\!\ }c@{\!\ }c@{\!\ }c@{\!\ }c@{\!}
c@{\!}c@{\!\ }c@{\!\ }c@{\!\ }c|}\hline
\ \ \ \ &&&&&&&&&&\ \ \\
&&&&&&&&\ \ &\C{\ast}&\\[-.3em]
&&&&&0&\cdots&0&&&\\[-.5em]
&&&&&\vdots&&\vdots&&&\\[-.2em]
&&&&{\ast}&0&&&&&\\[-.5em]
&&0&\cdots&0&&&\vdots&&&\\[-.5em]
&&\vdots&&&&&0&&&\\[-.0em]
&0 &0&&&&&\C{1}&&&\\[-0em]
&0&&&&&&0&&&\\[-.5em]
&\vdots&&&0&&&\vdots&&&\\[-.2em]
&0&&&&&&0&&&\\[-.2em]
&a\!\!+\!\!1\ &0&\cdots&&\cdots&0&b\!\!-\!\!1&&&\\[1em]
\hline\end{array}
\end{array}
$$
\\
Here the $\C{1}$ on the left is at position $(i_0,j_0)$,
\, $a\sh=m_{i_1j_0}$, and $b\sh=m_{i_1j_1}\sh>0$.
\\[1em]
(iii)(b)  The transpose of move (iii)(a).
Suppose $(i_0,j_0)\sh<(i_1,j_1)$
with  $(i_0,j_0)\sh\in\D$,\, $m_{i_0j_0}\sh=1$,\ $m_{i_1j_1}\sh>0$,
and $m_{ij}\sh=0$ whenever:
$(i_0,j_0)\sh<(i,j)\sh<(i_1,j_1)$,\, $(i,j)\sh\neq(i_0,j_1)$; 
and whenever $(i,j)\sh\leq(i_1,j_0)$, $(i,j)\!\not\leq\!\D$.
Then define:
$$
M':=M-E_{i_0j_0}-E_{i_1j_1}+E_{i_0j_1}+E_{i_1j_0},\qquad 
\D':=[\,\D\cup\{(i_1,j_0)\}\,]\,.
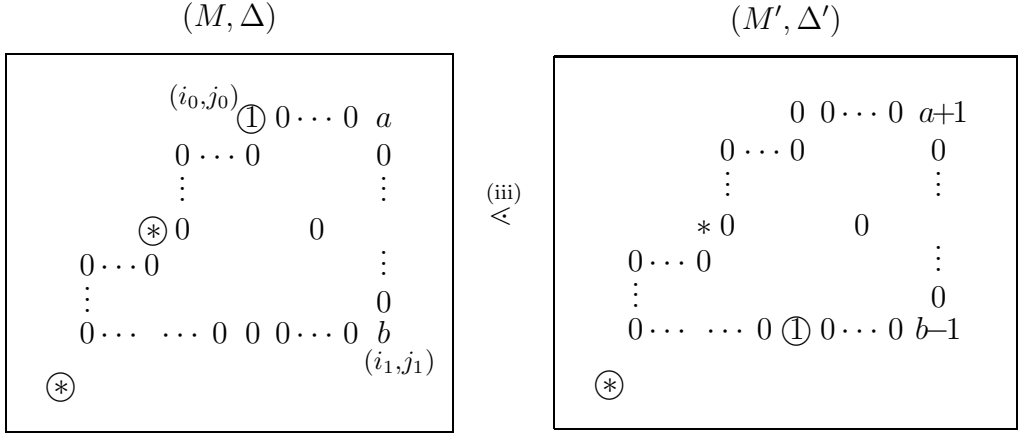
$$
$$\hspace{-.2in}
\begin{array}{c}(M,\D)\\[.5em]
\begin{array}{|c@{\!\ }c@{\!\ }c@{\!\ }c@{\!}c@{\!}c@{\!}c@{\!}
c@{\!\ }c@{\!\ }c@{\!\ }c@{\!\ }c@{\!\ }c@{\!\ }c@{\!\ }c|}\hline
\ \ &&&&&&&&&&&&&&\ \ \\[-.5em]
&&&&&&\hspace{-1em}\mbox{\footnotesize $(i_0,\!j_0)$}&&&&&&&&\\[-.5em]
&&&&&&&\C{1}\,&0&\cdots&0&a&&&\\[-.0em]
&&&&&0&\cdots&0&&&&0&&&\\[-.5em]
&&&&&\vdots&&&&&&\vdots&&&\\
&&&&\C{\ast}&0&&&&0&&&&&\\[-.5em]
&&0&\cdots&0&&&&&&&\vdots&&&\\[-.5em]
&&\vdots&&&&&&&&&0&&&\\[-.2em]
&&0&\cdots&&\cdots&0&0&0&\cdots&0&b&&&\\[-.4em]
&&&&&&&&&&&\mbox{\footnotesize $(i_1,\!j_1)$}\hspace{-1em}&&&\\[-.4em]
&\C{\ast}&&&&&&&&&&&&&\\[.8em]
\hline\end{array}
\end{array}
\ \move{(iii)}\
\begin{array}{c}(M',\D')\\[.5em]
\begin{array}{|c@{\!\ }c@{\!\ }c@{\!\ }c@{\!\ }c@{\!}c@{\!}c@{\!}
c@{\!\ }c@{\!\ }c@{\!\ }c@{\!\ }c@{\!\ }c@{\!\ }c@{\!\ }c|}\hline
\ \ &&&&&&&&&&&&&&\ \ \\
&&&&&&&0&0&\cdots&0&\ a\!\!+\!\!1&&&\\[-.0em]
&&&&&0&\cdots&0&&&&0&&&\\[-.5em]
&&&&&\vdots&&&&&&\vdots&&&\\
&&&&{\ast}&0&&&&0&&&&&\\[-.5em]
&&0&\cdots&0&&&&&&&\vdots&&&\\[-.5em]
&&\vdots&&&&&&&&&0&&&\\[-.2em]
&&0&\cdots&&\cdots&0&\C{1}\,&0&\cdots&0&b\!\!-\!\!1&&&\\[.5em]
&\C{\ast}&&&&&&&&&&&&&\\[.8em]
\hline\end{array}
\end{array}
$$ \\
Here the $\C{1}$ on the left is at position $(i_0,j_0)$,
\, $a\sh=m_{i_0j_1}$, and $b\sh=m_{i_1j_1}\sh>0$.
\\[1em]
(iv)(a) Suppose $i_0\sh<i_2\sh<i_1$ and $j_2\sh<j_0\sh<j_1$
with $(i_0,j_0),(i_2,j_2)\sh\in\D$,\, $m_{i_0j_0}\sh=m_{i_2j_2}\sh=1$,\, 
$m_{i_1j_1}\sh>0$, and $m_{ij}\sh=0$ whenever:  
$(i_0,j_2)\sh<(i,j)\sh<(i_1,j_1)$ and
$(i,j)\!\not\leq\!\D$ and $(i,j)\sh\neq(i_0,j_1),(i_1,j_2)$.
Then define:
$$
M':=M-E_{i_0j_0}-E_{i_1j_1}-E_{i_2j_2}+
E_{i_1j_2}+E_{i_2j_0}+E_{i_0j_1},\quad
\D':=[\, \D\cup\{(i_2,j_0)\}\,]\,.
$$
$$
\hspace{-.6in}
\begin{array}{c}(M,\D)\\[.5em]
\begin{array}{|c@{\!\ }c@{\!\ }c@{\!\ }c@{\!\ }c@{\!\ }c@{\!\ }c
@{\!\ }c@{\!\ }c@{\!\ }c@{\!\ }c@{\!\ }c@{\!\ }c|}\hline
&&&&&&&&&&&&\\[-.5em]
\ \ &&&&&&\mbox{\footnotesize $(i_0,\!j_0)$}\!&&&&&&\ \ \\[-.5em]
&&&&&&&\C{1}\,&0&\cdots&0&\ b&\\
&&&&&0&\cdots&0&&&&0&\\[-.5em]
&&&&&\vdots&&&&&&\vdots&\\[-.2em]
&&&&\C{\ast}\,&0&&&&&&&\\
&&0&\cdots&0&&&&&&&&\\[-.5em]
&&\vdots&&&&&&&&&&\\[-.2em]
\mbox{\footnotesize $(i_2,\!j_2)$}&\C{1}\,&0&&&&&0&&&&&\\
&0&&&&&&&&&&&\\[-.5em]
&\vdots&&&&&&&&&&\vdots&\\[-.2em]
&0&&&&&&&&&&0&\\[-.2em]
&a&0&\cdots&&&&&&\cdots&0&c&\\[-.8em]
&&&&&&&&&&&&\!\mbox{\footnotesize $(i_1,\!j_1)$}\\[.5em]
\hline\end{array}
\end{array}
\ \move{(iv)}\ 
\begin{array}{c}(M',\D')\\[.5em]
\begin{array}{|c@{\!\ }c@{\!\ }c@{\!\ }c@{\!\ }c@{\!\ }c@{\!\ }c
@{\!\ }c@{\!\ }c@{\!\ }c@{\!\ }c@{\!\ }c@{\!\ }c|}\hline
\qquad&&&&&&&&&&&&\qquad \\[-.0em]
&&&&&&&0&0&\cdots&0&\ b\!\!+\!\!1&\\
&&&&&0&\cdots&0&&&&0&\\[-.5em]
&&&&&\vdots&&&&&&\vdots&\\[-.2em]
&&&&\ast\ &0&&&&&&&\\
&&0&\cdots&0&&&&&&&&\\[-.5em]
&&\vdots&&&&&&&&&&\\[-.2em]
&0&0&&&&&\C{1}&\mbox{\footnotesize $(i_2,\!j_0)$}\!\!\!\!\!\!\!&&&&\\
&0&&&&&&&&&&&\\[-.5em]
&\vdots&&&&&&&&&&\vdots&\\[-.2em]
&0&&&&&&&&&&0&\\[-.2em]
&\!\!a\!\!+\!\!1\,&0&\cdots&&&&&&\cdots&0&\ c\!\!-\!\!1&\\[-.8em]
&&&&&&&&&&&&\\[.5em]
\hline\end{array}
\end{array}
$$ \\
Here the coordinate markings indicate that on the left, the
$\C{1}$'s are at positions $(i_0,j_0)$, $(i_2,j_2)$,
and $c=m_{i_1j_1}\sh>0$.  The $\C{1}$ on the right
is at position $(i_2,j_0)$.
\\[1em]
(iv)(bc)  Suppose we have $(i_0,j_0)\sh<(i_1,j_1)$
with $m_{i_0j_0},m_{i_1,j_1}\sh>0$, and
$m_{ij}\sh=0$ whenever: $(i_0,j_0)<(i,j)\sh<(i_1,j_1)$,
except for $(i,j)\sh= (i_0,j_1),(i_1,j_0)$ and one other
position as specified below.
Further suppose one of the following cases:\\
(b) we have $i_0\sh<i_2\sh<i_1$ with $(i_2,j_0)\sh\in\D$,\,
$m_{i_2j_0}\sh=1$, and $(i_0,j_1)\!\not\leq\!\D$; or\\
(c) we have $j_0\sh<j_2\sh<j_1$ with $(i_0,j_2)\sh\in\D$,\,
$m_{i_0j_2}\sh=1$, and $(i_1,j_0)\!\not\leq\!\D$.\\
Then define:
$$
M':=M-E_{i_0j_0}-E_{i_1j_1}+E_{i_0j_1}+E_{i_1j_0},\qquad \D':=\D\,.
$$
This is the same as move (ii), except that it occurs in the presence
of a $\C{1}$ at $(i_0,j_2)$ or $(i_2,j_0)$.
\\[.3em]
(b)
$$
\begin{array}{c} (M,\D)\\[.5em]
\begin{array}{|c@{\!\ \ }c@{\!\ \ }c@{\!\ \ }c@{\!\ \ }c@{\!\ \ \ }c@{\!\ \ }c|}
\hline
\ \ &&&&&&\ \ \\
&\ a\ & 0&\cdots & 0&\ b\ &\\[-.3em]
&0&&&&0&\\[-.5em]
&\vdots&&&&\vdots&\\
&0&&&&&\\[-.3em]
&\C{1}&&&&&\\[-0em]
&0&&&&&\\[-.5em]
&\vdots&&&&\vdots&\\[-.3em]
&0&&&&0&\\[-.3em]
& c&0&\cdots&0& d&\\[1em]
\hline
\end{array}
\end{array}
\ \move{(iv)}\  
\begin{array}{c} (M',\D')\\[.5em]
\begin{array}{|c@{\!\ }c@{\!\ \ }c@{\!\ \ }c@{\!\ \ }c@{\!\ \ }c@{\!\ }c|}
\hline
\ \ &&&&&&\ \ \\
& a\!\!-\!\!1 & 0&\cdots & 0& b\!\!+\!\!1&\\[-.3em]
&0&&&&0&\\[-.5em]
&\vdots&&&&\vdots&\\
&0&&&&&\\[-.3em]
&\C{1}&&&&&\\[-0em]
&0&&&&&\\[-.5em]
&\vdots&&&&\vdots&\\[-.3em]
&0&&&&0&\\[-.3em]
& c\!\!+\!\!1&0&\cdots&0& d\!\!-\!\!1&\\[1em]
\hline
\end{array}
\end{array}
$$ \\
Here $b$ must not be circled, and there must be no circled element
weakly southeast of $b$.  The move does not change the circled
positions.\\[.3em]
(c)
$$
\begin{array}{c} (M,\D)\\[.5em]
\begin{array}{|c@{\!\ \ }c@{\!\ }c@{\!\ }c@{\!\ }c@{\!\ }c@{\!\ }
c@{\!\ }c@{\!\ }c@{\!\ }c@{\!\ \ }c|}
\hline
\ \ &&&&&&&&&&\ \ \\
&\ a\ & 0&\cdots&0&\C{1}\ &0&\cdots & 0&\ b\ &\\[-.3em]
&0&&&&&&&&0&\\[-.5em]
&\vdots&&&&&&&&\vdots&\\
&0&&&&&&&&0&\\[-.3em]
& c&0&\cdots&&&&\cdots&0& d&\\[1em]
\hline
\end{array}
\end{array}
\ \move{(iv)}\  
\begin{array}{c} (M',\D')\\[.5em]
\begin{array}{|c@{\!\ \ }c@{\!\ }c@{\!\ }c@{\!\ }c@{\!\ }c@{\!\ }
c@{\!\ }c@{\!\ }c@{\!\ \ }c@{\!\ \ }c|}
\hline
\ \ &&&&&&&&&&\ \ \\
&a\!\!-\!\!1& 0&\cdots&0&\C{1}\ &0&\cdots & 0&b\!\!+\!\!1&\\[-.3em]
&0&&&&&&&&0&\\[-.5em]
&\vdots&&&&&&&&\vdots&\\
&0&&&&&&&&0&\\[-.3em]
& c\!\!+\!\!1&0&\cdots&&&&\cdots&0& d\!\!-\!\!1&\\[1em]
\hline
\end{array}
\end{array}
$$
\\
Here $c$ must not be circled, and there must be no circled element
weakly southeast of $c$.  The move does not change the circled
positions.
\\[1em]
(v)  Suppose, for $t\sh\geq 1$, we have $i_0\sh<i_1\sh<\cdots\sh<i_t$, and 
$j_1\sh>j_2\sh>\ldots\sh>j_t\sh>j_0$, with 
$(i_1,j_1),\ldots,(i_t,j_t)\sh\in\D$,\, $m_{i_0j_0}\sh>0$, and $m_{ij}\sh=0$ 
whenever: $(i_0,j_0)<(i,j)\leq(i_s\sh-1,j_s\sh-1)$ for some $s\sh=1,\ldots,t$. Then define: 
$$ 
M':=M-E_{i_0j_0}-E_{i_1j_1}-\cdots-E_{i_tj_t}
+E_{i_0j_1}+E_{i_tj_0}\,,
$$
$$
\D':=D\,
\setminus\{(i_0,j_0),(i_1,j_1),\ldots,(i_t,j_t)\}
\,\cup\{(i_0,j_1),(i_t,j_0)\}\,.
$$
\\[-1em]
$$
\hspace{-1.25in}
\begin{array}{c}(M,\D)\\[.5em]
\begin{array}{|c@{\!\ }c@{\!\ }c@{\!\ }c@{\!\ }c@{\!\ }c@{\!\ }c@{\!\ }c
@{\!\ }c@{\!\ }c@{\!\ }c@{\!\ }c@{\!\ }c@{\!\ }c@{\!\ }c@{\!\ }c|}\hline
\quad&&&&&&&&&&&&&&&\quad\\[.1em]
\mbox{\footnotesize $(i_0,\!j_0)$}\!\!\!\!&&&&&&&&&&&&&&&\\[-.6em]
&\ a\ &0&\cdots&&&&&&&&&\cdots&0&b&\\[-.3em]
&0&&&&&&&&&&&&&*&\\[-.5em]
&\vdots&&&&&&&&&&&&&\vdots&\\[-.3em]
&&&&&&&&&&&&&&*&\\[-.2em]
&&&&&0&&&&&&c&\cdots&*&\C{g}&\\[-.5em]
&&&&&&&&&&&&&&&\!\mbox{\footnotesize $(i_1,\!j_1)$}\\[-1.3em]
&&&&&&&&&&&\vdots&&&&\\[-.3em]
&&&&&&&&&&&*&&&&\\[-.2em]
&&&&&&&&d&\cdots&*&\C{h}&&&&\\[-.5em]
&&&&&&&&&&&&\!\mbox{\footnotesize $(i_2,\!j_2)$}&&&\\[-1.3em]
&&&&&&&&\vdots&&&&&&&\\[-.3em]
&&&&&&&&*&&&&&&&\\[-.2em]
&&&&&e&\cdots&*&\C{\mbox{\small \,$i$}}&&&&&&&\\[-.5em]
&&&&&&&&&\!\mbox{\footnotesize $(i_{t-1},\!j_{t-1})$}\!\!\!\!\!\!\!\!&&&&&&\\[-1.3em]
&\vdots&&&&\vdots&&&&&&&&&&\\[-.3em]
&0&&&&*&&&&&&&&&&\\[-.2em]
&f&*&\cdots&*&\C{\mbox{\small $j$}}&&&&&&&&&&\\[-.5em]
&&&&&&\!\mbox{\footnotesize $(i_t,\!j_t)$}&&&&&&&&&\\[.6em]
\hline\end{array}
\end{array}
\move{(v)}
\begin{array}{c}(M',\D')\\[.5em]
\begin{array}{|c@{\!\ }c@{\!\ }c@{\!\ }c@{\!\ }c@{\!\ }c@{\!\ }c@{\!\ }c
@{\!\ }c@{\!\ }c@{\!\ }c@{\!\ }c@{\!\ }c@{\!\ }c@{\!\ }c@{\!\ }c|}\hline
&&&&&&&&&&&&&&&\\[-.5em]
\!\!\!\!&&&&&&&&&&&&&&\mbox{\footnotesize $(i_0,\!j_1)$}&\\[-.0em]
&a\sh-1\,&0&\cdots&&&&&&&&&\cdots&0&\,b\sh+\C{1}&\\[-.2em]
&0&&&&&&&&&&&&&*&\\[-.5em]
&\vdots&&&&&&&&&&&&&\vdots&\\[-.3em]
&&&&&&&&&&&&&&*&\\[-1.3em]
&&&&&&&&&&\!\!\!\!\!\!\!\mbox{\footnotesize $(i_1,\!j_2)$}\!\!\!&&&&&\\[-0.3em]
&&&&&0&&&&&&c\sh+1&\cdots&*&{g\sh-1}&\\[-.4em]
&&&&&&&&&&&\vdots&&&&\\[-.8em]
&&&&&&&&&&&*&&&&\\[-.2em]
&&&&&&&&d\sh+{1}&\cdots&*&{h\sh-1}&&&&\\[-.5em]
&&&&&&&&\vdots&&&&&&&\\[-.3em]
&&&&&&&&*&&&&&&&\\[-1.3em]
&&&&\!\!\!\!\!\!\!\!\!\!\!\mbox{\footnotesize $(i_{t-1},\!j_t)$}\!\!\!&&&&&&&&&&&\\[-0.3em]
&&&&&e\sh+1&\cdots&*&{\mbox{\small $\,i\sh-1$}}&&&&&&&\\[-.4em]
&\vdots&&&&\vdots&&&&&&&&&&\\[-.3em]
&0&&&&*&&&&&&&&&&\\[-.3em]
&f\sh+\C{1}\,&*&\cdots&*&{\mbox{\small $j\sh-1$}}&&&&&&&&&&\\[-.2em]
&\!\mbox{\footnotesize $(i_t,\!j_0)$}&&&&&&&&&&&&&&\\[.6em]
\hline\end{array}
\end{array}
$$
\\[1em]
Here $m\sh+\C{1}$ means $m\sh+1$ circled, and
$*$ means a value $m_{ij}\geq 0$, which
is unchanged by the move. 
Note that, in contrast to moves (i)--(iv), we have $\D'<\D$.

\subsection{Strategy of proof}

We will prove the Rank and Move Theorems
for triple flags by means of three ``chain lemmas.''
\begin{lemma}
\label{Schain1}\qquad\qquad
$(M,\D)\lmv (M',\D')\leadsto (M,\D)\ldeg (M',\D')$
\end{lemma}
For $(M,\D)\move{mv} (M',\D')'$, we will give an explicit 
degeneration of $(M',\D')$ to $(M,\D)$.
\begin{lemma}
\label{Schain2}\qquad\qquad
$(M,\D)\ldeg (M',\D')'\leadsto (M,\D)\lrk (M',\D')$ 
\end{lemma}
This will follow from general principles of algebraic geometry. 
\begin{lemma} 
\label{Schain3}\qquad\qquad
$(M,\D)\lrk (M',\D')\leadsto (M,\D)\lmv (M',\D')$ 
\end{lemma}
This is a purely combinatorial result.  Given $(M,\D)\lrk (M',\D')$,
we construct $(\tM,\tD)$ with $(M,\D)\move{mv} (\tM\,\tD)\lrk (M',\D')$.

\section{Proof of Theorem \ref{Abruhat}}\label{sec:abruhat2}
\label{Abruhat proof}

In order to prepare and illuminate the proofs of
Lemmas \ref{Schain1}--\ref{Schain3} for triple flags, 
we give the precisely analogous arguments for the 
classical case of two flags, thereby proving 
Theorem \ref{Abruhat}.

\begin{lemma}
\label{Achain1}\qquad\qquad 
$M\lmv M'\leadsto M\ldeg M'$
\end{lemma}

\noindent{\it Proof.}
Given $M\move{mv} M'$, 
it suffices to find a one-parameter algebraic
family of double flags $(B_\bdot(\tau),C_\bdot(\tau))$,
indexed by $\tau\in\kk$, such that:
$$\begin{array}{c}
(B_\bdot(\tau),C_\bdot(\tau))\in F_{M'}
\ \ \text{for}\ \ \tau\sh\neq 0\,,
\\[.5em]
(B_\bdot(0),C_\bdot(0))\in F_{M}\ .
\end{array}
$$
Consider a basis of $V$ indexed according to $M=(m_{ij})$ as:
$$
V=\langle\, e_{ijk}\mid (i,j)\in[1,q]\sh\times[1,r],
1\sh\leq k\sh\leq m_{ij}\,\rangle\,,
$$
and define a set of vectors indexed according to $M'=(m'_{ij})$,
$$
\{\,v_{ijk}(\tau)\mid (i,j)\in[1,q]\sh\times[1,r],
1\sh\leq k\sh\leq m'_{ij}\,\}\,
$$
as follows.  
Let us use the symbol $e_{ij\max}$ to mean that the third subscript
in $e_{ijk}$ has as large a value as possible, 
namely $\max\sh=m_{ij}$ ; and similarly $v_{ij\max}$ 
means $\max\sh=m'_{ij}$ .  Now let:
$$\begin{array}{rcl}
v_{i_0j_1\max}(\tau)&:=&e_{i_0j_0\max}+\tau e_{i_1j_1\max}\\[.2em]
v_{i_1j_0\max}(\tau)&:=&e_{i_0j_0\max}\\[.2em]
v_{ijk}(\tau)&:=&e_{ijk}\qquad \text{otherwise}.
\end{array}
$$
For $\tau\sh\neq 0$, let 
$B_i(\tau):=\langle v_{i'jk}\mid i'\leq i\,\rangle$,\
$C_j(\tau):=\langle v_{ij'k}\mid j'\leq j\,\rangle$.
For $\tau\sh\neq 0$, the set $\{v_{ijk}\}$ forms a basis of $V$,
and with respect to this basis $(B_\bdot(\tau),C_\bdot(\tau))\in F_{M'}$ .

Now define $B_i(0)\!:=\!\lim_{\tau\to 0} B_i(\tau)$ and
$C_j(0)\!:=\!\lim_{\tau\to 0} C_j(\tau)$,
the limits in the Zariski topology.\footnote{
Note that these limits are guaranteed to exist
by the properness of $\PP^{n-1}\sh\times
\Flag(\bb)\sh\times\Flag(\cc)$.
We do not need this general fact, however, since
we explicitly identify the limit $(B_\bdot(0),
C_\bdot(0))$.
}
We proceed to evaluate these limits, showing that, 
with respect to the basis $\langle e_{ijk}\rangle$, 
we have $(B_\bdot(0),C_\bdot(0))\in F_M$ .

As $\tau\sh\to 0$, the vectors $\{\, v_{i'jk}(\tau)\mid i'\sh\leq i_0\,\}$
remain linearly independent, and $v_{i_0j_1\max}(\tau)\to e_{i_0j_0\max}$,
so we have:
$$
B_{i_0}(0)=\langle\, 
v_{i'jk}(0)\mid i'\sh\leq i_0\,\rangle
=\langle\, e_{i'jk}\mid i'\sh\leq i_0\,\rangle\,.
$$
Furthermore, we can take
linear combinations of basis vectors to find,
for arbitrary $\tau\sh\neq 0$:
$$\begin{array}{rcl}
B_{i_1}(\tau)&=&
\langle e_{i_0j_0\max}+\tau e_{i_1j_1\max},\,e_{i_0j_0\max}\rangle
\\[.2em]
&&\oplus\ \ \langle\, v_{i'jk}\mid i'\leq i_1,\,
(i,j,k)\sh\neq(i_0,j_1,m'_{i_0j_1}),(i_1,j_0,m'_{i_1,j_0})\,\rangle 
\\[.5em]
&=&
\langle e_{i_0j_0\max},\,e_{i_1j_1\max}\rangle
\\[.2em]
&&\oplus\ \ \langle\, e_{i'jk}\mid i'\leq i_0,\,
(i,j,k)\sh\neq(i_0,j_0,m_{i_0j_0}),(i_1,j_1,m_{i_1,j_1})\,\rangle 
\\[.5em]
&=&\langle e_{i'jk}\mid i'\leq i_1\rangle\,.
\end{array}$$
Since the final basis is constant with respect
to $\tau$,
the limit space $B_{i_1}(0)$ exists and
has the same basis.
Similarly for any $i,j$ we have $B_i(0)=\langle e_{i'jk}\mid i'\sh\leq i\rangle$, 
$C_{j}(0)=\langle e_{ij'k}\mid j\sh\leq j'\rangle$.    \pfover 

\begin{lemma}
\label{Achain2}\qquad\qquad 
$M\ldeg M'\leadsto M\lrk M'$
\end{lemma}

\noindent{\it Proof.}  Since $M\ldeg M'$, there exists
an algebraic family with
$(B_\bdot(\tau),C_\bdot(\tau))\sh\in F_{M'}$
for $\tau\sh\neq 0$, and $(B_\bdot(0),C_\bdot(0))\in F_{M}$.

For any fixed $d$, the condition $\dim(B_i\cap C_j)\geq d$
defines a closed subvariety of $\Flag(\bb)\sh\times\Flag(\cc)$
(cut out by the vanishing of certain determinants in the homogeneous
coordinates of the flag varieties).
Hence if the generic elements of the family
satisfy $r_{ij}(M')\geq d$, then so does
the limit element at $\tau=0$.
That is, the rank numbers $r_{ij}(M')$ can 
only get larger upon degeneration.
\pfover

\begin{lemma}
\label{Achain3}\qquad\qquad $M\lrk M'\leadsto M\lmv M'$\,.\\
In fact, if $M\slrk M'$, then we can perform a simple 
move on $M$ to obtain a matrix $\tM$ satisfying 
$$
M\move{mv} \tM \lrk M'\,.
$$
\end{lemma}

\noindent
{\it Proof.}  Denote
$M=(m_{ij})$,\, $r_{ij}=r_{ij}(M)$, \,
$M'=(m'_{ij})$,\, $r'_{ij}=r_{ij}(M')$, so that
$r_{ij}\geq r'_{ij}$ by assumption.
Consider the lexicographically first position $(k_0,l_0)$ where the 
matrices differ: that is, $m_{ij}=m'_{ij}$ for $i\sh<k_0$ and
for $i\sh=k_0$, $j\sh<l_0$; and $m_{k_0l_0}>m'_{k_0l_0}$, making
$r_{k_0l_0}>r'_{k_0l_0}$.  Consider a rectangle
$[k_0,k_1]\sh\times [l_0,l_1]$ as large as possible
such that $r_{ij}>r'_{ij}$ for 
$(i,j)\in[k_0,k_1\sh-1]\sh\times[l_0,l_1\sh-1]$.
(Such a rectangle is not necessarily unique.)

Claim: $m_{i_1j_1}>0$ for some $(i_1,j_1)\in
[k_0+1,k_1]\sh\times[l_0+1,l_1]$.
Otherwise we would have $m_{ij}=0$ for all 
$(i,j)\in [k_0\sh+1,k_1]\sh\times[l_0\sh+1,l_1]$.
By the maximality of the rectangle, we have
$r_{il_1}=r'_{il_1}$ and $r_{k_1j}=r'_{k_1j}$
for some $(i,j)\sh\in [k_0,k_1\sh-1]\sh\times[j_0,j_2\sh-1]$.
By the definition of the rank numbers and the vanishing of the $m_{kl}$,
we have:
$$
r_{k_1l_1}=r_{il_1}+r_{k_1j}-r_{ij}
<r'_{il_1}+r'_{k_1j}-r'_{ij}
\leq r'_{k_1l_1}.
$$
This contradiction proves our claim.

We may assume the $(i_1,j_1)$ found above is 
$\leq$-minimal, so that $m_{ij}=0$ for 
$(i,j)\in[k_0\sh+1,i_1]\sh\times[l_0\sh+1,j_1]$,
$(i,j)\sh\neq(i_1,j_1)$.
By moving right or down from 
$(k_0,l_0)$ to a position $(i_0,j_0)$, 
we can get:
$$
m_{i_0j_0},m_{i_1j_1}>0,\quad
m_{ij}=0 \ \ \text{for}\ \ 
(i_0,j_0)\sh<(i,j)\sh<(i_1,j_1), (i,j)\sh\neq(i_0,j_1),(i_1,j_0)
$$
\vspace{-1.3em}
$$
r_{ij}>r'_{ij}\ \ \text{for}\ \
(i,j)\sh\in[i_0,i_1\sh-1]\sh\times[j_0,j_1\sh-1]\,.
$$
These are all the necessary conditions to
perform our simple move: 
$\tM:=M-E_{i_0j_0}-E_{i_1j_1}+E_{i_1j_0}+E_{i_0j_1}$. 
Denoting $\tr_{ij}=r_{ij}(\tM)$, we have:
$$
\tr_{ij} =
\left\{\begin{array}{cl}
r_{ij}\sh-1\geq r'_{ij}&\text{if}\ (i,j)\in 
  [i_0,j_1\sh-1]\sh\times[j_0,j_1\sh-1]\\
r_{ij}\geq r'_{ij}&\text{otherwise,}
\end{array}\right.
$$
showing that $\tM\lrk M'$.  \pfover
\\

Thus, all three orders are identical, denoted $M\leq M'$.
It only remains to prove that the simple moves give
the covers of this order.

Suppose $M\move{R} M'$ for a rectangle 
$R=[i_0,i_1]\sh\times[j_0,j_1]$, and
$M\sh< \tM\sh\leq M'$.  Without loss of 
generality, we may assume $M\move{S} \tM$ for some
other rectangle $S$.
Then it is clear that 
$$
r_{ij}(M)\geq r_{ij}(\tM)\geq r_{ij}(M')\,.
$$
But 
$$
r_{ij}(M)=\left\{\begin{array}{cl}
r_{ij}(M')\sh+1 &\text{for}\ (i,j)\sh\in[i_0,i_1\sh-1]\sh\times
[j_0,j_1\sh-1]\\
r_{ij}(M')&\text{otherwise}\,,
\end{array}\right.
$$
and similarly for $S$.  Hence $S\subset R$, but because
$m_{ij}=0$ for all $(i,j)\in R$, we must have $S=R$,
and $\tM=M'$.

This completes the proof of Theorem \ref{Abruhat}

\section{Geometry of decorated matrices}\label{geom}

In this section,
we prove the first two chain lemmas 
for triple flags, the ones involving geometric arguments.
We follow the same lines of argument as for two flags.

Let $(M,\D)$, $(M',\D')$ be decorated 
matrices indexing
the orbits in $\PP^{n-1}\sh\times
\Flag(\bb)\sh\times\Flag(\cc)$, 
and denote as above
$M=(m_{ij})$, $r_{ij}=r_{ij}(M)$, etc., as well as
$\rb{ij}=r_{\langle ij\rangle}(M,\D)$, etc.
(Later we will also use $\tr_{ij}=r_{ij}(\tM)$, etc.)
Let us first prove the Uncircling Lemma:
\\[.5em]
{\it Proof of Lemma \ref{uncirc}.}\ Suppose we have
$(i_0,j_0),(i_1,j_1)\in S$ with $(i_0,j_0)<(i_1,j_1)$.  
Let $S'= S\setminus \{(i_0,j_0)\}$, and denote by 
$F_{M,S}\,, F_{M,S'}$ the $\GL(V)$-orbits 
defined analogously to $F_{M,\D}$.  It suffices to show that 
$F_{M,S}=F_{M,S'}$ .

Let $\langle e_{ijk}\rangle$ be a basis of $V$ indexed
according to $M$ (that is, $k\leq m_{ij}$),
and let $B_i:=\langle e_{i'jk}\mid i'\leq i\rangle$,
$C_j:=\langle e_{ij'k}\mid j'\leq j\rangle$,
and $A:=\langle\, \sum_{(i,j)\in S} e_{ij1}\,\rangle$.
By definition, $(A,B_\bdot,C_\bdot)\in F_{M,S}$ .

Now let us write the same $(A,B_\bdot,C_\bdot)$ in terms of
a new basis $\langle v_{ijk}\rangle$, defined by
$v_{i_1j_11}:=e_{i_1j_11}+e_{i_0j_01}$ 
and $v_{ijk}=e_{ijk}$ otherwise.
We have 
$$\begin{array}{rcl}
B_{i_1}&=&\langle\, e_{i_0j_01},\, e_{i_1j_11}\,\rangle\\[.2em]
&&\oplus\ \langle\, e_{ijk}\mid i\leq i_1,\, (i,j,k)\sh\neq
(i_0,j_0,1),(i_1,j_1,1)\,\rangle\\[.5em]
&=&\langle\, e_{i_0j_01},\, e_{i_1j_11}+e_{i_0j_01}\,\rangle\\[.2em]
&&\oplus\ \langle\, e_{ijk}\mid i\leq i_1,\, (i,j,k)\sh\neq
(i_0,j_0,1),(i_1,j_1,1)\,\rangle\\[.5em]
&=&\langle\, v_{ijk}\mid i\leq i_1\,\rangle
\end{array}$$
Similarly $B_i=\langle v_{i'jk}\mid i'\leq i\rangle$ and
$C_j=\langle v_{ij'k}\mid j'\leq j\rangle$ for all $i,j$.
Finally,  
$$\begin{array}{rcl}
A&=&\langle\, e_{i_0j_01}+e_{i_1j_1}
+\sum_{(i,j)\in S''} e_{ij1}\,\rangle\\[.5em]
&=&\langle\, v_{i_1j_1}+\sum_{(i,j)\in S''} v_{ij1}\,\rangle
\ \ =\ \ \langle\, \sum_{(i,j)\in S'} v_{ij1}\,\rangle\,,
\end{array}$$
where $S''=S\setminus\{(i_0,j_0),(i_1,j_1)\}$.
Hence $(A,B_\bdot,C_\bdot)\in F_{M,S'}$ , proving the Lemma. 
\pfover
\\[1em]
{\it Proof of Lemma \ref{Schain1}.}\ \  We follow the method
in the proof of Lemma \ref{Achain1}.  For each type of move 
$(M,\D)\move{mv}(M',\D')$ listed in Sec.~\ref{moves}, we start with a basis 
$V=\langle\,e_{ijk}\mid k\leq m_{ij}\,\rangle$
indexed according to $M$;
then define another set of vectors 
$\langle\,v_{ijk}(\tau) \mid k\leq m'_{ij}\,\rangle$
indexed according to $M'$
which forms a basis when $\tau\sh\neq 0$.

We then define a family of triple flags for 
$\tau\sh\neq 0$ by:
$$
B_i(\tau):=\langle v_{i'jk}(\tau)\mid i'\leq i\rangle,\ \ 
C_j(\tau):=\langle v_{ij'k}(\tau)\mid j'\leq j\rangle\,,
$$
$$
A(\tau):=\left\langle\, \sum_{(i,j)\in \D'} v_{ij1}(\tau)\,\right\rangle\,,
$$
so that by definition $(A(\tau),B_\bdot(\tau),C_\bdot(\tau))
\in F_{M',\D'}$ with respect to the basis 
$\langle v_{ijk}(\tau)\rangle$. 
Furthermore, it will be clear that as $\tau\to 0$, we have 
$$
(A(\tau),B_\bdot(\tau),C_\bdot(\tau))\to
(A(0),B_\bdot(0),C_\bdot(0))\,,
$$
where 
$$
B_i(0):=\langle e_{i'jk}\mid i'\leq i\rangle,\ \ 
C_j(0):=\langle e_{ij'k}\mid j'\leq j\rangle\,.
$$
$$
A(0):=\left\langle\, \sum_{(i,j)\in \D} e_{ij1}\,\right\rangle\,.
$$
This shows that $\lim_{\tau\to 0}(A(\tau),B_\bdot(\tau),C_\bdot(\tau))\in F_{M,\D}$ ,
and proves the Lemma.  We will give details only in the last, most
complicated case.
\\[.5em]
(i)  Let $S:=\{(i,j)\in\D\mid (i,j)<(i_1,j_1)\}$.  Now define:
$v_{i_1j_11}:=\tau e_{i_1j_11}+\sum_{(i,j)\in S}e_{ij1}$ and 
$v_{ijk}:=e_{ijk}$ otherwise.  (In each of the following cases,
we implicitly define $v_{ijk}:=e_{ijk}$ for those $(i,j,k)$ which
are not otherwise specified.)
\\[.5em]
(ii) Since $\D=\D'$, the line $A$ is unchanged in the move.
Thus we may take $v_{ijk}$ exactly 
as in the proof of Lemma \ref{Achain1}.
\\[.5em]
(iii) Let $S:=\{(i,j)\in\D\mid (i,j)<(i_0,j_1)\}$.  For (iii)(a),
define: $v_{i_1j_0\max}:=e_{i_0j_01}$
and $v_{i_0j_11}:=\tau e_{i_1j_1\max}+\sum_{(i,j)\in S}e_{ij1}$ .
Transpose for (iii)(b).
\\[.5em]
(iv)(a) Let $S:=\{(i,j)\in\D\mid (i,j)<(i_2,j_0)\}$. 
Now define:  $v_{i_1j_2\max}:=e_{i_2j_21}+\tau e_{i_1j_1\max}$,\ \
$v_{i_0j_1\max}:=e_{i_0j_01}+\tau e_{i_1j_1\max}$, and
$v_{i_2j_01}:=\sum_{(i,j)\in S}e_{ij1}$.
\\[.5em]
(iv)(bc)  Same as (ii).
\\[.5em]
(v)  Let $S:=\D\setminus\{(i_1,j_1),\ldots,(i_t,j_t)\}$.
Define:
$$v_{ij1}:=\tau e_{ij1}\ \ \text{for}\ (i,j)\,\in\,S\, ;\qquad
v_{i_sj_sk}:=e_{i_s,j_s,k+1}\ \ \text{for}\ s=1,2,\ldots, t\,;
$$
$$
v_{i_0j_11}:= e_{i_0j_0\max}
+\tau\,\mbox{$\sum_{s=1}^t e_{i_sj_s1}$}\,, \qquad
v_{i_0j_1k}:=e_{i_0,j_1,k-1}\ \ \text{for}\ k>1\,;
$$
$$
v_{i_tj_01}:= -e_{i_0j_0\max}\,, \qquad
v_{i_tj_0k}:=e_{i_t,j_0,k-1}\ \ \text{for}\ k>1\,;
$$
$$
\begin{array}{r@{\!}l}
v_{i_{s}j_{s+1}\max}:=e_{i_0j_0\max}&\,+\,\tau^2 e_{i_1j_11}
+\tau^2 e_{i_2j_21}+\cdots+\tau^2 e_{i_sj_s1}\\
&+\tau e_{i_{s+1}j_{s+1}1}+\cdots+\tau e_{i_tj_t1}
\end{array}
\ \ \text{for}\ s=1,2,\ldots,t\sh-1\,. 
$$
The crucial part of the transition matrix 
between $\langle v_{ijk}\rangle$ and
$\langle e_{ijk}\rangle$, containing all the nonzero
``non-diagonal'' coefficients, is:
$$\begin{array}{c|cccccc}
&v_{i_0j_11}&v_{i_1j_2\max}& v_{i_2j_3\max}&\cdots&v_{i_{t-1}j_t\max}&v_{i_tj_01}\\
\hline\\[-.9em]
 e_{i_0j_0\max}&1&1&1&\cdots&1&\!\!-1\ \\
e_{i_1j_11}&\tau&\tau^2&\tau^2&\cdots&\tau^2&0\\
e_{i_2j_21}&\tau&\tau&\tau^2&\cdots&\tau^2&0\\
\vdots&\vdots&\vdots&\vdots&&\vdots&\vdots\\
e_{i_{t-1}j_{t-1}1}&\tau&\tau&\tau&\cdots&\tau^2&0\\
e_{i_tj_t1}&\tau&\tau&\tau&\cdots&\tau&0
\end{array}
$$
Let us focus on these basis vectors, suppressing 
all other vectors with ellipsis marks \ldots.
We compute, as $\tau\to 0$ :
$$
B_{i_0}(\tau)=\langle v_{i_0j_11},\ldots\rangle
\to\langle e_{i_0j_0\max},\ldots\rangle\,,
$$
$$\begin{array}{rcl}
B_{i_1}(\tau)&=&\langle v_{i_0j_11},v_{i_1j_2\max},\ldots\rangle
=\langle v_{i_0j_11},\,\tau^{-1}(v_{i_0j_11}-v_{i_1j_2\max}),\ldots\rangle\\[.2em]
&=&\langle v_{i_0j_11},\,(1-\tau)e_{i_1j_11},\ldots\rangle
\to\langle e_{i_0j_0\max},e_{i_1j_11},\ldots\rangle\,,
\end{array}$$
and similarly for the other $B_i$.  Also: $$ 
C_{j_0}(\tau)=\langle 
v_{i_tj_01},\ldots\rangle \to\langle e_{i_0j_0\max},\ldots\rangle \,,
$$ 
$$\begin{array}{rcl} C_{j_t}(\tau)&=&\langle 
v_{i_tj_01},v_{i_{t-1}j_t\max},\ldots\rangle =
\langle -v_{i_tj_01},\,
\tau^{-1}(v_{i_{t-1}j_t\max}+v_{i_tj_01}),\ldots\rangle\\[.2em] 
&=&\langle v_{i_tj_01},\,\tau e_{i_1j_11}+\cdots
+\tau e_{i_{t-1}j_{t-1}1}+e_{i_tj_t1},\ldots\rangle \to
\langle e_{i_0j_0\max},e_{i_tj_t1},\ldots\rangle \,,
\end{array}$$
and similarly for the other $C_j$.  Finally, 
$$\begin{array}{rcl}
A(\tau)&=&\langle\, v_{i_0j_11}+v_{i_tj_01}
+\sum_{(i,j)\in S}v_{ij1}\,\rangle\\[.4em]
&=&\langle\, \tau e_{i_1j_1}+\cdots+\tau e_{i_tj_t1}
+\sum_{(i,j)\in S}\tau e_{ij1}\,\rangle\\[.4em]
&=&\langle\, e_{i_1j_1}+\cdots+e_{i_tj_t1}
+\sum_{(i,j)\in S}e_{ij1}\,\rangle\,.
\end{array}$$
In each case the $\tau\to 0$ limit is the desired subspace,
an element of $F_{M,\D}$ 
with respect to the basis $\langle e_{ijk}\rangle$. 
\pfover
\\[1em]
{\it Proof of Lemma \ref{Schain2}.}\
We must show that $\rb{ij}(A,B_\bdot,C_\bdot)\geq d$ is a closed
condition on $\PP(V)\sh\times\Flag(\bb)\sh\times\Flag(\cc)$.
This is clear if $i\sh=0$ or $j\sh=0$.
Given $(i,j)\in [1,q]\sh\times[1,r]$ and
a triple flag $X=(A,B_\bdot,C_\bdot)$, define the linear map
$\phi_{ij}^X:B_i\!\times C_j\to V/A$\ \ by 
$(v_1,v_2)\mapsto v_1+v_2\mod A$.
Then $\phi_{ij}^X$ depends algebraically on $X$:
indeed, we may write $\phi_{ij}^X$ in coordinates as a matrix
of size $(n\sh-1)\sh\times (b_i\sh+c_j)$ with entries depending
polynomially on the homogeneous coordinates of $X$.
Thus $\dim\Ker(\phi_{ij}^X)\geq d$ is a closed condition on $X$.
However, 
$$\begin{array}{rcll}
\dim\Ker(\phi_{ij}^X)&=& \dim B_i+\dim  C_j-\mathop{\rm rank}(\phi_{ij}^X)\\[.3em]
&=& \dim B_i+\dim C_j-\dim(\,(B_i\sh+C_j)/A\,)\\[.3em]
&=& \dim B_i+\dim C_j-\dim(B_i\sh+C_j)+\dim(\,A\cap(B_i\sh+C_j)\,)\\[.3em]
&=& \dim(B_i\cap C_j)+\dim(\,A\cap(B_i\sh+C_j)\,)\\[.3em]
&=&\rb{ij}(A,B_\bdot,C_\bdot)&\mbox{\pfover}
\end{array}$$

\section{Combinatorics of decorated matrices}\label{combin}

To prove the Rank and Move Theorems, it only remains to show
Lemma \ref{Schain3}.  

Thus, assume we are given
$(M,\D)\slrk (M',\D')$.  We wish to show that
we can perform one of the 
moves (i)--(v) on $(M,\D)$ to obtain 
a decorated matrix $(\tM,\tD)$ satisfying:
$$
(M,\D)\move{mv} (\tM,\tD) \lrk (M',\D')\,.
$$
We give an algorithm producing $(\tM,\tD)$.

For cases (i)--(v)(a), we assume
$\D\leq \D'$.  This is equivalent to 
$\delta_{ij}(\D)\geq
\delta_{ij}(\D')$ for all $(i,j)$.
Note that if $(\tM,\tD)$ is 
any decorated matrix satisfying $\tM\lrk M'$ and $\tD\leq \D'$, 
then $(\tM,\tD)\lrk(M',\D')$.
\\[1em]
{\sc case}\! (i)\
Assume $M=M'$.  Then $\D<\D'$, and we may choose a minimal
$(i_1,j_1)$ with $(i_1,j_1)\not\leq\D$, $(i_1,j_1)\leq\D'$.  
Now we may apply move (i) to the appropriate block 
$[i_0,i_1]\sh\times[j_0,j_1]$ of $(M,\D)$, obtaining $(\tM,\tD)$.
Since $M=\tM=M'$ and $\D<\tD\leq\D'$, we have: 
$$
(M,\D)\move{(i)}(\tM,\tD)\lrk(M',\D')\,.
$$
\\[-.5em]
\noindent
For the rest of the cases, 
we assume $M\neq M'$, so that $M\slrk M'$ 
in the rank order for two flags.  
We may then apply the proof 
of Lemma \ref{Achain3} to find 
$(i_0,j_0),(i_1,j_1)$ such that:
$$
m_{i_0j_0},m_{i_1j_1}>0,\quad m_{ij}=0\ \ \text{for}\ \
(i_0,i_1)\sh<(i,j)\sh<(i_1,j_1), (i,j)\sh\neq(i_0,j_1),(i_1,j_0)
$$
\vspace{-1.5em}
$$
\text{and}\quad r_{ij}>r'_{ij}\ \text{for}\ (i,j)\sh\in R,
\ \ \text{where}\ \ R:=[i_0,i_1\sh-1]\sh\times[j_0,j_1\sh-1]\,.
$$
Henceforth, in cases (ii)--(v)(a), 
we will assume as given such positions
$(i_0,j_0)$, $(i_1,j_1)$.
\\[1em]
{\sc case}\! (ii)\
Assume $(i_0,j_0),(i_1,j_1)\not\in\D$;
or $(i_0,j_0)\in\D$, $m_{i_0j_0}>1$.

If $(i_0,j_1),(i_1,j_0)\in\D$, then
we may apply move (i) to the block
$[i_0,i_1]\sh\times[j_0,j_1]$ of $(M,\Delta)$,
obtaining $\tM=M$,\, $\tD=[\D\cup\{(i_1,j_1)\}]$.
Clearly for all $(i,j)$,\ 
$\tr_{ij}=r_{ij}\geq r'_{ij}$;
and for $(i,j)\not\in $,\ 
$\trb{ij}=\rb{ij}\geq\rb{ij}'$; 
while for $(i,j)\in R$, 
$\trb{ij}=r_{ij}\geq r_{ij}'\sh+1\geq\rb{ij}'$.
Thus $(M,\D)\move{(i)}(\tM,\tD)
\leq(M',\D').$

If, on the other hand, $(i_0,j_1)\not\in\D$
or $(i_1,j_0)\not\in\D$, then we may
apply move (ii) to the block
$[i_0,i_1]\sh\times[j_0,j_1]$ of $(M,\Delta)$,
obtaining $(M,\D)\move{(ii)}(\tM,\tD)$.  
For any $(i,j)$ we have $\tM\lrk M'$ by Lemma \ref{Achain3},
and by assumption $\tD=\D\leq\D'$, so 
$(M,\D)\move{(ii)}(\tM,\tD)\lrk(M',\D')$.
\\[1em]
{\sc case}\! (iii)(a)\  Assume $(i_0,j_0)\in\D$, $m_{i_0j_0}=1$, 
and $(i_0,j_1)\leq\D'$.

If $m_{ij}>0$, for some $(i,j)$ in the rectangle
$[1,i_0]\sh\times[j_0,j_1]$ and with
$(i,j)\not\leq\D$, we may apply move (i) to some
block in this rectangle, obtaining 
$(M,\D)\move{(i)}(\tM,\tD)$
with $M=\tM\lrk M'$ and $\D<\tD\leq\D'$.  
Thus $(M,\D)\move{(i)}(\tM,\tD)\lrk(M',\D')$.  

Otherwise $m_{ij}=0$ for all $(i,j)$ in
$[1,i_0]\sh\times[j_0,j_1]$ with $(i,j)\not\leq\D$, 
so we have all the conditions necessary to
apply move (iii)(a) to the block
$[i_0,i_1]\sh\times[j_0,j_1]$ of $(M,\Delta)$,
obtaining $(M,\D)\move{(iii)}(\tM,\tD)$.  
As before, $\tM\lrk M'$ by Lemma \ref{Achain3},
and $\tD\leq\D\sh\cup\{(i_0,j_1)\}\leq\D'$ by assumption,
so $(M,\D)\move{(iii)}(\tM,\tD)\lrk(M',\D')$.
\\[1em]
{\sc case}\! (iii)(b)\  Assume $(i_0,j_0)\in\D$, $m_{i_0j_0}=1$, 
and $(i_1,j_0)\leq\D'$.  We can repeat the arguments 
of the previous case with rows and columns transposed.
\\[1em]
{\sc case}\! (iv)(a)\  Assume $(i_0,j_0)\in\D$, $m_{i_0j_0}=1$,
and $(i_1,j_0),(i_0,j_1)\not\leq\D'$; and furthermore
that $m_{ij}=m'_{ij}$ for $i\sh\leq i_0$ and for $i\sh=i_0,j\sh<j_0$.
(In terms of the proof of Lemma \ref{Achain3}, the last condition means 
$(i_0,j_0)=(k_0,l_0)$.)

Our assumptions imply:
$$
\D'\leq\{(i_0\sh-1,r),(q,j_0\sh-1),(i_1\sh-1,j_1\sh-1)\}.
$$
Since $(i_0,j_0)\in\D\leq\D'$, there must exist
$(s_1,t_1)\in\D'\sh\cap R$.

Claim: $r_{i,j_0-1}>r_{i,j_0-1}'$ for $i\in[s_1,i_1\sh-1]$.  
Otherwise, if $r_{i,j_0-1}=r_{i,j_0-1}'$ 
for some $i\in[s_1,i_1\sh-1]$, 
we would have:
\begin{eqnarray*}
r_{it_1}&=& r_{i,j_0-1}+r_{i_0-1,t_1}-r_{i_0-1,j_0-1}+1\\
&\leq& r'_{i,j_0-1}+r'_{i_0-1,t_1}-r'_{i_0-1,j_0-1}+
\!\!\!\!\sum_{(i_0,j_0)\leq(k,l)\leq(i,t_1)}\!\!\!\!\!\!\!\! m'_{kl} \\[-.8em]
& = & r'_{it_1}\,.
\end{eqnarray*}
This contradicts $r_{ij}>r'_{ij}$ within $R$, establishing the claim.

Now take a rectangle $S:=[s_0,i_1\sh-1]\sh\times[t_0,j_0\sh-1]$ as large
as possible such that $r_{ij}>r'_{ij}$ for $(i,j)\in S$.  
Further, take $(s_0,t_0)$ to be lexicographically minimal with
the above property, so that $s_0\leq s_1$ by the above Claim.

Claim: There exists $(i_2,j_2)\in S$
such that $m_{i_2j_2}>0$.  Otherwise, we would
have $m_{ij}=0$ within $S$.  By the maximality of $S$,
there exists $(i,j)\in S$ so that $r_{i,t_0-1}=r'_{i,t_0-1}$
and $r_{s_0-1,j}=r'_{s_0-1,,j}$.  Then
$$
r_{ij}=r_{s_0-1,j}+r_{i,t_0-1}-r_{s_0-1,t_0-1}
\leq r'_{s_0-1,j}+r'_{i,t_0-1}-r'_{s_0-1,t_0-1}
\leq r'_{ij}\,.
$$
This contradicts $r_{ij}>r'_{ij}$ within $S$, establishing the claim.

Thus, we may choose $(i_2,j_2)$ in $S$ with $m_{i_2j_2}>0$
and $m_{ij}=0$ for $(i_2,j_2)<(i,j)\leq(i_1\sh-1,j_0\sh-1)$.
In fact, choose $(i_2,j_2)$ to be as northeast as possible
with these properties, so that $m_{ij}=0$ for 
$(i,j)\in[s_0,i_2\sh-1]\sh\times[j_2\sh+1,j_0\sh-1]$.
If $(i_2,j_2)\not\in\D$ or $m_{i_2j_2}>1$,
we have all the necessary conditions 
to apply move (ii) 
to the block $[i_2,i_1]\sh\times[j_2,j_1]$ or
$[i_2,i_1]\sh\times[j_2,j_0]$ 
of $(M,\D$), and finish as
in case (ii) above, obtaining $$(M,\D)\move{(ii)}
(\tM,\tD)\lrk(M',\D')\,.$$

If $m_{ij}>0$ for some $(i,j)$ in the rectangle
$[1,s_0\sh-1]\sh\times[j_2,j_0\sh-1]$ with $(i,j)\not\leq\D$, 
then we may apply move (i) to some block within this rectangle.
Since the southeast corner $(s_0\sh-1,j_0\sh-1)\leq(s_1,t_1)\in\D'$,
we obtain:  $$(M,\D)\smallstack{(i)}{<}
(\tM,\tD)\lrk(M',\D')\,.$$

Thus, we have reduced to the case where
$(i_2,j_2)\in\D$, $m_{i_2j_2}=1$, and $m_{ij}=0$ 
or $(i,j)\in\D$ for all 
$(i,j)\in[1,i_2\sh-1]\sh\times[j_2,j_0\sh-1]$.

If $m_{i_1j_0}>0$, we may apply move (iii)
to the block $[i_2,i_1]\sh\times[j_2,j_0]$ of $(M,\D$).  
Then clearly $\tr_{ij}\geq r'_{ij}$, since 
the block lies inside $S$.
For $(i,j)$ outside $[i_0,i_2\sh-1]\sh\times[j_2,j_0\sh-1]$, 
we have 
$$
(i\sh+1,j\sh+1)\leq\tD\Leftrightarrow
(i\sh+1,j\sh+1)\leq\D\Rightarrow(i\sh+1,j\sh+1)\leq\D'\,,
$$
so 
$$
\trb{ij}=\tr_{ij}+\delta_{ij}(\tD)=
\tr_{ij}+\delta_{ij}(\D) 
\geq r'_{ij}+\delta_{ij}(\D')
=\rb{ij}'\,.
$$
For $(i,j)\in[i_0,s_0\sh-1]\sh\times[j_2,j_0\sh-1]$, 
we have 
$$
(i\sh+1,j\sh+1)\leq(s_0,j_0)\leq(s_1,t_1)\in\D'\,,
$$
and $\tr_{ij}=r_{ij}$,
so $\trb{ij}\geq r_{ij}\geq r'_{ij}=\rb{ij}'$.
For $(i,j)\in[s_0,i_2\sh-1]
\sh\times[j_2,j_0\sh-1]\subset S$, 
we have $\trb{ij}\geq\tr_{ij}=r_{ij}\geq r'_{ij}+1\geq\rb{ij}'$.
Therefore $$(M,\D)\smallstack{(iii)}{<}
(\tM,\tD)\lrk(M',\D')\,.$$

Finally, suppose $m_{i_1j_0}=0$
or $(i,j)\in\D$,
and as before $(i_2,j_2)\in\D$,\, $m_{i_2j_2}=1$,
and $m_{ij}=0$ for all $(i,j)\in[i_0\sh+1,i_2\sh-1]\sh\times[j_2\sh+1,j_0\sh-1]$.
Then we may apply move (iv)(a) to the rectangle 
$[i_2,i_1]\sh\times[j_2,j_1]$ in $(M,\D)$.  Now, 
for $(i,j)$ outside the region
$$
[i_2,i_1\sh-1]\sh\times[j_2,j_1\sh-1]\ 
\cup\ [i_0,i_1\sh-1]\sh\times[j_0,j_1\sh-1]\subset S\cup R
$$
we have $\tr_{ij}=r_{ij}\geq r'_{ij}$; whereas 
for $(i,j)$ inside this region we have $\tr_{ij}=r_{ij}-1\geq r'_{ij}$.

To check $\trb{ij}\geq \rb{ij}'$, we repeat the argument
we used to show $(M,\D)\smallstack{(iii)}{<}$ $
(\tM,\tD)\lrk(M',\D')$ immediately above.
For $(i,j)$ outside the rectangle $[i_0,i_2\sh-1]\sh\times[j_2,j_0\sh-1]$,
we have $(i\sh+1,j\sh+1)\leq\D\Leftrightarrow (i\sh+1,j\sh+1)\leq\tD
\Rightarrow(i\sh+1,j\sh+1)\leq\D'$, so clearly $\trb{ij}\geq\rb{ij}'$
as above.  Similarly repeat the arguments for $(i,j)$ inside 
$[i_0,i_2\sh-1]\sh\times[j_2,j_0\sh-1]$, to obtain
 $$(M,\D)\smallstack{(iv)}{<}
(\tM,\tD)\lrk(M',\D')\,.$$
{\sc case}\! (iv)(bc)\ Assume $(i_0,j_0)\in\D$, $m_{i_0j_0}=1$,
and $(i_1,j_0),(i_0,j_1)\not\leq\D'$; but
that $m_{ij}\neq m'_{ij}$ for some $i\sh\leq i_0$ or for $i\sh=i_0,j\sh<j_0$.
From the proof of Lemma \ref{Achain3}, the last condition means 
that there exists $k_0<i_0$ with $m_{k_0j_0}>0$, or 
there exists $l_0<j_0$ with $m_{i_0l_0}>1$.  
Assume the first alternative 
(the other one being merely a transpose). 

By increasing $k_0$ if necessary, 
we may assume that
$m_{ij}=0$ for $(i,j)\in [k_0,i_1]\sh\times[j_0,j_1]$
except that definitely $m_{ij}>0$ for 
$(i,j)=(k_0,j_0),(i_0,j_0),(i_1,j_1)$,
and possibly $m_{ij}>0$ for
$(i,j)=(i_1,j_0),(i_0,j_1)$ and for $i=k_0$.

Suppose $m_{i_0j_1}>0$. If $m_{k_0j}>0$ for $j_0<j<j_1$;
or if $m_{k_0j_0}>0$ and $(k_0,j_1)\not\in\D$;
then we may apply move (ii) to the rectangle $[k_0,i_0]\sh\times[j,j_1]$
$[k_0,i_0]\sh\times[j_0,j_1]$, 
and finish as in case (ii):
$(M,\D)\move{(ii)}(\tM,\tD)\lrk(M',\D')$.
Otherwise, if $m_{k_0j_0}>0$,\
$m_{k_0j}=0$ for $j_0<j<j_1$,\ and
$(k_0,j_1)\in\Delta$; then apply
move (i) to $[k_0,i_0]\sh\times[j_0,j_1]$
to get $\tD=[\D\cup\{(i_0,j_1)\}]$, yielding
$(M,\D)\move{(i)}(\tM,\tD)\lrk(M',\D')$.

Now suppose $m_{i_0j_1}=0$.
If $m_{k_0j}>0$ for $j_0<j<j_1$;
then we may apply move (ii) to the rectangle $[k_0,i_1]\sh\times[j,j_1]$:
$(M,\D)\move{(ii)}(\tM,\tD)\lrk(M',\D')$.
If $m_{k_0j_1}\in\D$,
then we may apply move (iii) to 
$[i_0,i_1]\sh\times[j_0,j_1]$, getting
$(M,\D)\move{(iii)}(\tM,\tD)\lrk(M',\D')$.

Suppose none of the above cases hold.
Then $m_{i_0j_1}=0$,\ $m_{k_0j_1}\not\in\D$, 
and $m_{k_0j}=0$ for $j_0<j<j_1$.
Then we finally have the conditions to
apply move (iv)(b) to the rectangle
$[k_0,j_0]\sh\times[j_0,j_1]$, and finish by again 
noting that $\tM\lrk M'$, $\tD=\D\leq\D'$,
so that $(M,\D)\move{(iv)}(\tM,\tD)\lrk(M',\D')$.
\\[1em]
{\sc case}\! (v)(a)\  Assume $(i_1,j_1)\in\D$.  
Then perform move (v) on 
the rectangle $[i_0,i_1]\sh\times[j_0,j_1]$ to obtain
$(M,\D)\smallstack{(v)}{<}(\tM,\tD)$.  We clearly have
$\tM\lrk M'$ as well as $\tD\leq\D\leq\D'$.
Hence $(\tM,\tD)\lrk(M',\D')$.
\\[1em]
We have now proved the Lemma assuming $\D\leq\D'$.
Next assume $\D\not\leq\D'$.
\\[.5em]
{\sc case}\! (v)(b)\  Consider any $(i,j)$ with $(i\sh+1,j\sh+1)\leq\D$ but
$(i\sh+1,j\sh+1)\not\leq\D'$.  Then we have
$r_{ij}+0=\rb{ij}\geq\rb{ij}'=r'_{ij}+1$, so that $r_{ij}>r'_{ij}$.

 Now, choose some $(k_1,l_1)\in\D$ with $(k_1,l_1)\not\leq\D'$; and
choose $(k_0,l_0)<(k_1,l_1)$ so that the rectangle $[k_0,k_1]\sh\times
[l_0,l_1]$ is as large as possible with $r_{ij}>r'_{ij}$ for
$(i,j)\in[k_0,k_1\sh-1]\sh\times[l_0,l_1\sh-1]$.
(The above remarks show that $k_0 <k_1$,\,
$l_0 <l_1$.)

Claim:  There exists $(i_0,j_0)\in
[k_0,k_1\sh-1]\sh\times[l_0,l_1\sh-1]$ 
with $m_{i_0j_0}>0$.  Otherwise, we would have $m_{ij}=0$
inside this rectangle.  By the maximality, 
there exists $(i,j)$ in the rectangle with
$r_{i,l_0-1}=r'_{i,l_0-1}$ and $r_{k_0-1,j}=r'_{k_0-1,j}$.
Then 
$$
r_{k_0-1,l_0-1}=r_{k_0-1,j}+r_{i,l_0-1}-r_{ij}
<r'_{k_0-1,j}+r'_{i,l_0-1}-r'_{ij}
\leq r'_{k_0-1,l_0-1}\,.
$$
This contradicts $r_{ij}\geq r'_{ij}$, 
establishing the claim.

Let us say a rectangle $[i,k]\sh\times[j,l]$
has property (R) if:
$$
(\mbox{\rm R})\qquad 
\begin{array}{c}
m_{ij},m_{kl}>0,\qquad
(k,l)\in\D,\qquad
(k,l)\not\leq\D',
\\[.5em]
r_{st}>r'_{ij}\ \ \text{for}\ \ (s,t)\sh\in
[i,k\sh-1]\sh\times[j,l\sh-1],
\end{array}$$
Our argument above shows that 
we may choose $(i_0,j_0)$ so that 
$[i_0,k_1]\sh\times[j_0,l_1]$
has property (R).
Next choose points $(i_1,j_1),
(i_2,j_2),\ldots,(i_t,j_t)$,
as many as possible, such that: 
the rectangles $[i_0,i_s]\sh\times[j_0,j_s]$
all have property (R), and
$$
\begin{array}{c}
i_1<i_2<\cdots<i_t\,,\qquad
j_1>j_2>\cdots>j_t\,,
\\[.3em]
[i_s,i_{s+1}\sh-1]\sh\times[j_{s+1},j_s\sh-1]
\ \cap\ \D =\emptyset
\\[.3em]
[i_s,i_{s+1}\sh-1]\sh\times[j_{s+1},j_s\sh-1]
\ \cap\ \D' =\emptyset
\end{array}
$$
for all $s$.  That is, $(i_1,j_1),(i_2,j_2),\ldots,(i_t,j_t)$ are 
consecutive elements of $\D$ 
(listed from northeast to southwest), 
not separated by intermediate elements of $\D'$.
By moving $(i_0,j_0)$ southeast if necessary,
and keeping only the rectangles with
$(i_0,j_0)\leq(i_s,j_s)$,
we may assume, in addition to all the above,
that $m_{ij}=0$ for all
$(i,j)\neq(i_0,j_0)$ in the region
$$
S:=\bigcup_{s=1}^t\ [i_0,i_s\sh-
1]\sh\times[j_0,j_s\sh-1]\,,
$$
the union of all our rectangles.   

Now apply move (v) to the region
$[i_0,i_1]\sh\times[j_0,j_t]$ of
$(M,\D)$, which contains all the
$[i_0,i_s]\sh\times[j_0,j_s]$ as subrectangles, 
and obtain $(M,\D)\smallstack{(v)}{<}
(\tM,\tD)$.  
Then
$\tr_{ij}=r_{ij}$ outside $S$,
and $\tr_{ij}=r_{ij}-1\geq r'_{ij}$ inside $S$.
Hence $\tr_{ij}\geq r'_{ij}$ for all $(i,j)$.

Further, for $(i,j)\not\in S$, we have
$\delta_{ij}(\tD)=\delta_{ij}(\D)$,
so clearly $\trb{ij}=\rb{ij}\geq \rb{ij}'$.

Finally, consider $(i,j)\in S$.
The definition of $\tD$ and $S$ ensures that 
$\delta_{ij}(\tD)=0\Rightarrow
\delta_{ij}(\D')=0$;
if this were not so, we could enlarge the list 
$(i_1,j_1),\ldots,(i_t,j_t)$ while 
keeping $(i_0,j_0)$ fixed.
Hence for $(i,j)\in S$, we have:
$$
\trb{ij}=\tr_{ij}+\delta_{ij}(\tD)\geq
r'_{ij}+\delta_{ij}(\D')=\rb{ij}'\,.
$$
That is, $\trb{ij}\geq\rb{ij}'$ 
for all $(i,j)$, and $(\tM,\tD)\lrk(M',\D')$.

This concludes the proof of Lemma \ref{Schain3}, and hence
of the Rank and Move Theorems.

\section{Minimality}

In this section, we prove 
the Minimality Theorem for triple flags.
In the case of full flags $\bb=\cc=(1^n)$
considered in the Introduction,
this follows from the fact that each simple
move corresponds to a codimension-one
containment of orbits 
$F_{M,\D}\subset\overline F_{M',\D'}$.
However, this is not true in general.
We give an alternative purely combinatorial argument.

By the Move Theorem, only the moves
(i)--(v) are candidates for covers of
our Bruhat order:
we must show that each such move is indeed
a minimal relation.

We will denote the Bruhat order on decorated matrices
simply by $(M,\D)\leq(M',\D')$.
We retain the notations
$m_{ij}, m'_{ij},\tm_{ij}$ for the 
matrix entries of $M, M',\tM$; 
and $r_{ij}, r'_{ij}, \tr_{ij}$ for the rank-numbers
$r_{ij}(M),r_{ij}(M'),r_{ij}(\tM)$.

Given a decorated matrix $(M,\D)$,
suppose we perform some simple move on the 
block $[i_0,i_1]\sh\times[j_0,j_1]$, obtaining $(M',\D')$.
We say a position $(i,j)\in[1,q]\sh\times[1,r]$
is {\it $M$-inactive} with respect to the move
$(M,\D)<(M',\D')$ if: 
$$
r_{ij}\sh=r'_{ij},\ r_{i-1,j}\sh=r'_{i-1,j},\
r_{i,j-1}\sh=r'_{i,j-1},\ r_{i-1,j-1}=r'_{i-1,j-1}\,.
$$
Otherwise $(i,j)$ is {\it $M$-active}.
Similarly, $(i,j)$ is {\it $\D$-inactive} if:
$$
\rb{i-1,j-1}\sh=\rb{i-1,j-1}',\ 
\rb{i-1,j}\sh=\rb{i-1,j}',\
\rb{i,j-1}\sh=\rb{i,j-1}'\,;
$$
and otherwise $(i,j)$ is {\it $\D$-active}.
The following result implies that we may obtain $(M',\D')$ from $(M,\D)$
by changing the entries of $M$ only at the $M$-active positions,
and the elements of $\D$ only at the 
$M$-active and $\D$-active positions,
leaving $M$ and $\D$ unchanged at all inactive positions.

\begin{lemma} If $(i,j)$ is $M$-inactive, 
then $m_{ij}\sh=m'_{ij}$.
If $(i,j)$ is $M$-inactive and $\D$-inactive, then
$(i,j)\in\D \Leftrightarrow (i,j)\in\D'$.
\end{lemma}

\noindent{\it Proof.} For the first statement,
$m_{ij}=r_{ij}-r_{i-1,j}-r_{i,j-1}+r_{i-1,j-1}$.
For the second statement, note that $(i,j)\in\D$
exactly when $(i,j)\leq\D$ and $(i\sh+1,j),(i,j\sh+1)\not\leq\D$.
\pfover
\\[1em]
Now suppose we have a possible couterexample 
to the Minimality Theorem,
a single move relation $(M,\D)\move{mv}(M',\D')$,
with an intervening element $(\tM,\tD)$
in the Bruhat order: 
$$(M,\D)<(\tM,\tD)\leq(M',\D')\,.$$
We must show that $(\tM,\tD)=(M',\D')$.

Since rank numbers must decrease at an active
position, we have:  if $(i,j)$ is $M$-active for the move
$(M,\D)<(\tM,\tD)$, then $(i,j)$ is $M$-active for the move
$(M,\D)<(M',\D')$; and similarly for $\D$-active positions.
That is, the move $(M,\D)\move{mv}(\tM,\tD)$
may act only at the active positions of $(M,\D)<(M',\D')$.

We now need only inspect the active positions
of each possible move $(M,\D)\move{mv}(M',\D')$, 
and verify that
there is no other move 
$(M,\D)\move{mv}(\tM,\tD)$ which
can be performed on the 
given active positions 
and which satisfies $(\tM,\tD)\slrk(M',\D')$.  
\\[.3em]
(i) If $(M,\D)\smallstack{(i)}{<}(M',\D')$, then 
there are no $M$-active positions, so that $\tM=M$.
Hence $\D<\tD<\D'$, but this is clearly impossible
since $m_{ij}=0$ for all $(i,j)$ with
$(i,j)\not\leq\D$,\ $(i,j)\leq\D'$,\ $(i,j)\neq(i_1,j_1)$.
\\[.3em]
(ii) The active positions 
form the block $[i_0,i_1]\sh\times[j_0,j_1]$. 
If $\D$ is disjoint from this block, 
then only other moves of type (ii) are possible, but
there are none such.  

If $(i_0,j_0)\in\D$, $m_{i_0j_0}>1$, we must also
consider type (i) applied to $\{i_0\}\sh\times[j_0,j_1]$
or $[i_0,i_1]\sh\times\{j_0\}$ or $[i_1,i_2]\sh\times[j_0,j_1]$,
each of which leads to $(M,\D)\smallstack{(i)}{<}
(\tM,\tD)\not<(M',\D')$, since
$\trb{i_0-1,j_0}=r_{i_0-1,j_0}<r_{i_0-1,j_0}+1=\rb{i_0-1,j_0}'$.
Also possible is type (iii)(a) applied to 
$[i_0,i_1]\sh\times[j_0,j_1]$, 
but then $(\tM,\tD)\not<(M',\D')$,
as above.  We apply the transposed arguments for (iii)(b).

Similarly for the cases $(i_0,j_1)\in\D$ and $(i_1,j_0)\in\D$.
\\[.3em]
(iii)(a) The $M$-active positions are the block
$[i_0,i_1]\sh\times[j_0,j_1]$, and the
$\D$-active positions are those $(i,j)\in
[i_0,i_1]\sh\times[1,j_1]$ with $(i,j)\not\leq\D$.
The only other possible moves are of type (i).
If we apply move (i) to any block inside 
$[i_0,i_1]\sh\times[1,j_0\sh-1]$, then we must
have $\D<\tD<\D'$, which is impossible.
If we apply (i) to $[i_0,i_1]\sh\times[j_0,j_1]$, 
then $(M,\D)\smallstack{(i)}{<} 
(\tM,\tD)\not<(M',\D')$ as before.
Transpose for (iii)(b).  
\\[.3em]
(iv)(a)  The $M$-active positions are  
$
\bR=[i_0,i_1]\sh\times[j_0,j_1]\, \cup\, [i_2,i_1]\sh\times[j_2,j_1]\,,
$
whereas the $\D$-active positions are $(i,j)\in[i_2,i_1]\sh\times[j_0,j_1]$
with $(i,j)\not\leq\D$.
The other possible moves are (iii) applied to $[i_2,i_1]\sh\times
[j_0,j_2]$ or to $[i_2,i_0]\sh\times[j_0,j_1]$; or (i) applied
to these same blocks.   All of these give $(\tM,\tD)\not<(M',\D')$.
\\[.3em]
(iv)(b) The active positions are $[i_0,i_1]\sh\times[j_0,j_1]$.
One possible move is (i) applied to 
$[i_2,i_1]\sh\times[j_0,j_1]$.  In this case $\trb{i_2,j_0-1}
=\rb{i_2,j_0-1}\sh-1<\rb{i_2,j_0-1}=\rb{i_2,j_0-1}'$, 
so $(\tM,\tD)\not\leq(M',\D')$.

Another possibility is (iii)(a) applied to $[i_2,i_1]\sh\times[j_0,j_1]$.
Again $\trb{i_0-1,j_0}<\rb{i_0-1,j_0}'$ and $(\tM,\tD)\not\leq(M',\D')$.
Similarly for (iii)(b) applied to $[i_2,i_1]\sh\times[j_0,j_1]$.
\\[.3em]
(iv)(c) Transpose of (iv)(b)
\\[.3em]
(v)  The active positions are:
$
\bigcup_{s=1}^t [i_0,i_s]\sh\times[j_0,j_s]\,.
$
One possible move is (ii) applied to some
block $[i_0,i_s]\sh\times[j_0,j]$ for
$j_{s+1}\leq j\leq j_s$.
Then $\trb{i_0j_0}=\rb{i_0j_0}-1<\rb{i_0j_0}
=\rb{i_0j_0}'$,
so $(\tM,\tD)\not\leq(M',\D')$.
Similarly for (ii) applied to a block
$[i_0,i]\sh\times[j_0,j_s]$.

The only other possible move is (v) 
applied to some smaller
block  $[i_0,i_l]\sh\times[j_0,j_m]$, where $[l,m]\subset [1,t]$
(strict inclusion).  Then 
$\trb{i_{l-1}j_{m+1}}=\rb{i_{l-1}j_{m+1}}\sh-1\sh+1
<\rb{i_{l-1}j_{m+1}}\sh+1=\rb{i_{l-1}j_{m+1}}'$,
because $(i_0,j_0)<(i_{l-1},j_{m+1})$.
Hence $(\tM,\tD)\not\leq(M',\D')$.
\\

Minimality is thus proved.

\end{document}